%% file: stringy.revised.tex
\documentclass [12pt,oneside] {amsart}
\usepackage{latexsym,amsfonts,amssymb,euler,euscript,amscd,fullpage}
\usepackage[pdftex]{graphicx}

\input moonshine.macros.tex

\usepackage [all]{xy}
\SelectTips{cm}{12}
\title{Stringy power operations in Tate $K$-theory}
\author{Nora Ganter\\ The University of Melbourne and Colby College}
\thanks{The author was supported by NSF grant DMS-0504539. Most
  of this paper was written during her stay at MSRI in the spring of
  2006 and the following year at University of Illinois at Urbana
  Champaign.} 
\date {\today}
\renewcommand {\AA} {{\hat A}}
\newcommand {\ABS} {Atiyah-Bott-Shapiro }

\newcommand {\ageg} {{\on{age}(g)}}
\newcommand {\ageh} {{\on{age}(h)}}
\newcommand {\bfG} {{\bf G}}
\newcommand {\Ckg} {{C^k_g}}
\newcommand {\CRg} {{C_g^{\RR/k}}}
\newcommand {\Dev} {{\on{Dev}}}
\newcommand {\Ell} {{\on{Ell}}}
\newcommand {\eq} {{\on{eq}}}

\renewcommand {\gg} {{\langle g\rangle}}
\newcommand {\hkr} {Hopkins-Kuhn-Ravenel }
\newcommand {\inv} {^{-1}}

\newcommand {\kLM} {{{}_k\mathcal L(M)}}

\renewcommand {\LG} {{\mathcal {LG}}}
\newcommand {\Lgs} {{\mathcal L_{(\ug,\sigma)}M^n}}

\newcommand {\LMG} {{\mathcal L(M\mmod G)}}

\renewcommand {\LX} {{\mathcal L\underline X}}
\newcommand {\maps}{\on{maps}}
\newcommand {\mA} {{\mathcal A}}
\newcommand {\mG} {{\mathcal G}}
\newcommand {\mI} {{\mathcal I}}
\newcommand {\mL} {{\mathcal L}}
\newcommand {\mP} {{\mathcal P}}
\newcommand {\mV} {{\mathcal V}}

\newcommand {\MGSn} {{(M\mmod G)\wr\Sn}}
\newcommand {\MSn} {{M^n\mmod\Sn}}
\newcommand {\MUP} {{\on{MUP}}}

\newcommand {\qjl} {{q^\frac jl}}

\newcommand {\sdim} {{\on{sdim}}}
\renewcommand {\SO} {\on{SO}}
\newcommand {\Spin} {{\on{Spin}}}
\newcommand {\Spinc} {{\on{Spin^\CC}}}
\renewcommand {\SS} {{{\mathbb S}^1}}
\renewcommand {\Stop} {{\Sym^\top_t}}
\newcommand {\str} {{\on{string}}}
\newcommand {\symsn} {{\on{sym}_\str^n}}
\renewcommand {\tilde}[1] {\widetilde{#1}}
\newcommand {\ts} {{(\underline\tau,\sigma)}}
\newcommand {\tts} {{(\tau_1,\dots,\tau_n,\sigma)}}
\newcommand {\ug} {{\underline g}}
\newcommand {\uh} {{\underline h}}
\newcommand {\ui} {{\underline i}}
\newcommand {\uj} {{\underline j}}
\newcommand {\ux} {{\underline x}}

\newcommand {\LgM} {{\mL_gM}}

\newcommand {\Vrest} {{V_{g\neq 1}}}
\newcommand {\VV} {{V}}
\newcommand {\LV} {{\mathcal L\VV}}

\newcommand {\Rl} {{\RR/l\ZZ}}
\newcommand {\RR} {{\mathbb R}}

\renewcommand {\XX} {{\underline X}}
\newcommand {\zlj} {{\zeta_l^j}}
\begin{document}
\begin{abstract} 
We study the loop spaces of the symmetric powers of an orbifold and
use our results to
define equivariant power operations in Tate $K$-theory. We prove that
these power operations are elliptic and 
that the Witten genus is an $\hinf$-map. As a corollary, we recover
a formula by Dijkgraaf, Moore, Verlinde and Verlinde for the orbifold
Witten genus of these symmetric powers. We outline some of the
relationship between our power operations and notions from 
(generalized) Moonshine.
\end{abstract}
\maketitle
\tableofcontents
\section{Introduction}
In recent years, a string theoretic result by Dijkgraaf, Moore,
Verlinde, and Verlinde has found much attention by mathematicians. The
celebrated paper \cite{DMVV} considers the symmetric powers $\MSn$ of
a compact, closed, complex manifold $M$ and expresses an invariant of
these symmetric powers, called the {\em orbifold elliptic genus},
merely in terms of the 
elliptic genus of $M$ itself:
\begin{equation}\label{DMVV-Eqn}
  \sum_{n\geq 0}\forb\MSn t^n = 
  \exp\(\sum_{m\geq 1}T_m(\phi(M))t^m\).  
\end{equation}
Here the $T_m$ are Hecke operators. In \cite{Borisov:Libgober:elliptic},
\cite{Borisov:Libgober:McKay},  
Borisov and Libgober proved
\eqref{DMVV-Eqn}, using algebro-geometric methods, and proceeded to
prove a McKay correspondence result for $\phi_{\orb}$.
In \cite{Tamanoi:01}, Tamanoi made a connection to homotopy
theory. In \cite{Tamanoi:03} and
\cite{Tamanoi:Infinite}, he explored the geometric side of
Formula \eqref{DMVV-Eqn}, generalizing the result to equivariant
manifolds and higher genus world-sheets and giving a beautiful
mathematical account of the loop 
space picture 
underlying the geometry of \eqref{DMVV-Eqn}, which is closer to the original
argument in \cite{DMVV} and \cite{Dijkgraaf} than the proof in
\cite{Borisov:Libgober:DMVV}.  

In \cite{Ganter:thesis}, one finds another link to homotopy theory:
interpreting the elliptic genus as a natural transformation between
cohomology theories (namely complex cobordism and elliptic cohomology), 
I proved that Formula \eqref{DMVV-Eqn} holds
whenever the natural transformation preserves cohomology operations.

From this point of view it is no surprise that some of the non-equivariant
picture in \cite{DMVV} and \cite{Tamanoi:03}
was independently discovered by Lupercio, Uribe, and
Xicot\'encatl, in
\cite{Lupercio:Uribe:Xicotencatl:symmetric}, as part of a program to
define cohomology operations in Chas-Sullivan cohomology of (the Borel
construction of) orbifolds.

However striking the formal analogies between 
the homotopy theoretic notions introduced in \cite{Ganter:thesis} on
one hand and results from orbifold string theory on the other, in
\cite{Ganter:thesis}, they
remained analogies. The present paper aims to bridge between the
different points of view by studying a version of equivariant elliptic
cohomology introduced by Devoto, whose objective was to capture the
behavior of orbifold loop spaces.
%
%
%

To motivate our results, let us recall the setup of
\cite{Ganter:thesis} in more detail: 
let $\Ell$ be an elliptic cohomology theory as defined in
\cite{Ando:Hopkins:Strickland:cube}, and let
$$
  \{\Ell_G\mid G \text{ is a finite group}\}
$$
be a compatible family of equivariant versions of $\Ell$. Let
$\map\phi\MU\Ell$ be a map of ring spectra, and assume that we are
given equivariant versions $\map{\phi_G}{\MU_G}{\Ell_G}$ of $\phi$.
The reader not familiar with elliptic cohomology 
should merely note that, typically, one can interpret
elements $\chi\in\Ell_G(\pt)$ as class functions on pairs of commuting
elements of $G$, in which case we say that $\Ell$ has a \hkr
(HKR)
character theory. 
Two important examples of such elliptic cohomology theories with HKR-theory
are
Borel-equivariant Lubin-Tate-Morava $E_2$-theory \cite{Hopkins:Kuhn:Ravenel}
and Devoto's equivariant Tate $K$-theory \cite{Devoto}. 
The former was the original object of study in
\cite{Hopkins:Kuhn:Ravenel} and \cite{Ganter:thesis}; the latter will
be the framework for the paper at hand.
\begin{Def}[Orbifold elliptic genera]
  Assume that $\Ell_G$ has a \hkr theory, and let $\chi\in\Ell_G(\pt)$.
  We define
  $$\eps_G(\chi):=\GG G\sum_{gh=hg}\chi(g,h),$$
  and
  $$
    \phi_\orb := \eps_G\circ\phi_G.
  $$
\end{Def}
Often, $\eps_G$ and thus $\phi_\orb$ take values in the coefficient
ring $\Ell_*$. 
Consider the following diagram:
\begin{equation}
  \label{intro-comm-diagram-Eqn}
  \xymatrix{{\MU_*}\ar[0,2]^{\phi}\ar[d]_{\oplus P_n}& & 
  {\Ell_*}\ar@{..>}[d]^{\oplus P_n}\\
  {\bigoplus\MU_{\Sn *}t^n}\ar[0,2]_{\oplus\phi_\Sn} & & 
  {\bigoplus\Ell_{\Sn*}t^n}\ar@{..>}[0,2]_{\oplus \eps_\Sn t^n}&& {\Ell\ps t}
  }  
\end{equation}
Here the left vertical arrow is the total power operation in
cobordism,
$$
  P_n([M]) = [M^n\acted\Sn],
$$
and $t$ is a dummy variable. The horizontal dotted arrow exists if
$\Ell$ has a \hkr theory, and then the composite of the left vertical
arrow with the two lower horizontal arrows becomes the left-hand side
of the DMVV formula (\ref{DMVV-Eqn}).

Assume now that $\Ell$ possesses power operations
$$
  \map{P_n}{\Ell_G(X)}{\Ell_G(X^n).}
$$
For $X$ the one point space and $G$ the trivial group, this implies
that the (right) vertical dotted arrow of the above diagram exists. The map
$\phi$ is called an 
{\em $\hinf$-map} if for each $n$, the equality $$P_n^{\Ell}\circ\phi =
\phi_\Sn\circ P_n^\MU$$ holds. For such $\phi$, the
diagram commutes. 

If $\Ell$ has a \hkr theory, one can use the $P_n$ to define Hecke
operators
$$
  \map{T_n}{\Ell(X)}{\Ell(X)}
$$
(c.f.\ \cite{Ando:thesis}). The $P_n$ are called {\em elliptic} if
this definition of Hecke operators agrees with another one, which is
in terms of isogenies \cite{Ando:Hopkins:Strickland:sigma}. 
By \cite{Ganter:thesis}, the composite of the dotted arrows is
\begin{equation}
  \label{expTm-Eqn}
  \sum_{n\geq 0}\eps_\Sn\(P_n(x)\) t^n = \exp\(\sum_{m\geq1}T_m(x)t^m\).
\end{equation}
If $\phi$ is an $\hinf$-map, this implies formula
(\ref{DMVV-Eqn}).

\subsection{Plan}
In the present paper, we treat the case where $\Ell$ is
Devoto's equivariant Tate $K$-theory and $\phi$ is 
the equivariant Witten genus. Section \ref{Thom-Sec} introduces the
equivariant Witten genus as a map of spectra.
Motivated by the work of Lupercio, Uribe and Xicot\'encatl
on orbifold loop spaces, we define Thom classes in Devoto $K$-theory
for $G$-equivariant
$U^2$-bundles, where $G$ is a finite group. We prove that the induced map of
spectra realizes the $G$-equivariant Witten genus defined in
\cite{FNLU}. We prove that the associated orbifold genus is
Morita-invariant and takes values in $$q^{-d/12}K_\Tate(pt) =
q^{-d/12}\ZZ\ps 
q,$$ where $d$ is the complex dimension of the manifold.

We proceed to study loop spaces of symmetric powers of orbifolds: If
$M$ is a manifold, then a loop in $\MSn$ is given by $n$ paths in $M$
and a permutation $\sigma$,
where for each $i$, the end point of the path $\gamma_i$ is the
starting point of the path $\gamma_{\sigma(i)}$.
For each $k$-cycle of $\sigma$, this produces one loop of length $k$
in $M$. 
The fact that the loop space of
$\MSn$ is made up from loops in $M$ is the key argument in \cite{DMVV}.
For the symmetric powers of a (global quotient) orbifold $M\mmod G$, a
similar picture is true, but the situation is complicated by the
additional action of $G$.

In Section \ref{loop-spaces-of-symmetric-Sec}, we prove a theorem
expressing the orbifold loop spaces of the symmetric powers of $M\mmod
G$ in terms of the orbifold loop space of $M\mmod G$.
This is 
a variation on Tamanoi's result, and we claim no originality.
However, in
order to introduce the notation for Section
\ref{power-operations-Sec}, we give a complete proof.
We also show that
the construction is compatible with iterated symmetric powers. 

In Section \ref{power-operations-Sec}, we use the analysis of
\ref{loop-spaces-of-symmetric-Sec} to define equivariant power
operations in Devoto $K$-theory. 
We compute a formula for the Hecke operators associated to our power
operations and conclude that the $P_n$ are elliptic. Interestingly,
this agrees exactly with the formula
for the twisted Hecke operators in generalized Moonshine found in
\cite{Ganter:twistedhecke}. If $G$ is the trivial group, we recover
the operators of \cite{Ando:poweroperations}. 

We proceed to prove that a variant of the
equivariant Witten genus is an $\hinf$-map. 
As a corollary, we obtain
the Dijkgraaf-Moore-Verlinde-Verlinde formula \ref{DMVV-Eqn} for our equivariant
Witten genus.

The total symmetric power $\Sym^\str_t$ is defined as the composite of
the two dotted arrows in Diagram (\ref{intro-comm-diagram-Eqn}). It is
computed by 
Formula (\ref{expTm-Eqn}). If $\phi(M)$ is the Witten genus of $M$, then
$\Sym^\str_t(\phi(M))$ becomes the right-hand side of the DMVV-formula. On the
other hand, we compute 
$\Sym_t^{\str}(V)$, where $V$ is a complex vector bundle, and find that it
equals Witten's exponential characteristic class
$$
  \Sym^\str_t(V) = \bigotimes_{k\geq 1} \Sym_{t^k}(V),
$$
where $\Sym_t$ is the total symmetric power in $K$-theory.
This sheds some light on the fact that the roles of $t$ and $q$ become
symmetric, when
the right-hand side of (\ref{DMVV-Eqn}) is written in its product form
$$\prod_{i,j}\(\frac1{1-q^it^j}\)^{c(ij)}.$$
Here the $c(i)$ are the coefficients of the Witten genus of $M$,
$$\phi(M) = \sum c(i)q^i.$$

To conclude this introduction, I will comment in some more detail on
how this paper compares to \cite{Ganter:thesis} and \cite{DMVV}.

In \cite{Ganter:thesis} 
I worked in a purely homotopy-theoretic setup and applied my results to 
the $\sigma$-orientation
of Morava $E_2$-theory, which was shown to be $\hinf$ in
\cite{Ando:Hopkins:Strickland:sigma}. 
Morava $E$-theory has no known geometric definition, and the fact that
I could reproduce a DMVV-type formula for the $\sigma$-orientation
suggests a deep and somewhat mysterious connection between
homotopy-theoretic notions such as the $K(n)$-local categories and
string-theoretic notions such as orbifold genera. 

The spectrum $K_\Tate$ is related to both of these:
It is an elliptic spectrum, and Devoto's
equivariant versions fit into the general framework of equivariant
elliptic cohomology and level structures on elliptic curves,
axiomatically formalized in 
\cite{Ginzburg:Kapranov:Vasserot}. 
Further, when restricted to $\MU\langle6\rangle$, the Witten genus
becomes the 
$\sigma$-orientation of $K_\Tate$ (cf.\
\cite{Ando:Hopkins:Strickland:cube}). 

On the other hand, as we have explained in \cite{Ando:French:Ganter}, 
the Witten genus is closely related to the genus 
considered by 
Borisov and Libgober. 
Moreover, Devoto's definition of $K_\Dev$ was inspired by orbifold
loop spaces, and so are our definitions of Thom classes and power
operations. 
%
\subsubsection{Acknowledgements} Many thanks go to Haynes Miller,
who generously shared his unpublished work with me. 
Parts of this
paper are closely following his exposition in \cite{Haynes:fpt} and
\cite{Haynes:loops}.
I would like to thank Nitu Kitchloo for many helpful suggestions on an
earlier version of this paper. 
It is a pleasure to thank Ralph Cohen, Matthew Ando, Haynes Miller,
Alex Ghitza, Michael Hopkins and Ernesto Lupercio for long and helpful
conversations. 
\subsection{Notation index}

\newcommand{\ite}[1]{\item[#1 \hfill]}

\begin{list}{}
   {\setlength{\labelwidth}{3cm} 
    \setlength{\labelsep}{0cm}
    \setlength{\leftmargin}{3.5cm}
    \setlength{\rightmargin}{0cm}
    \setlength{\itemsep}{0cm} 
    \setlength{\parsep}{0.05cm} 
    \setlength{\itemindent}{0cm}
    \setlength{\listparindent}{0cm}}
\ite{$\ageg$} {see Page \pageref{age-page}},
\ite{$C_g=C_G(g)$} {the centralizer of $g$ in $G$},
\ite{$C_g^k$} {the group $C_g\times (\ZZ/k|g|\ZZ) / \langle
  g,-k\rangle$, c.f.\ Definition \ref{extended-centralizer-def}},
\ite{$C_g^{\mathbb R/k\ZZ}$} {the group $C_g\times (\mathbb R/k|g|\ZZ) / \langle
  g,-k\rangle$, c.f.\ Definition \ref{extended-centralizer-def}},
\ite{$\det$} {determinant line bundle (top exterior power) of a
  complex line bundle} 
\ite{$D_n$} {c.f.\ Definitions \ref{D_n-Def} and \ref{D_n-Def2}},
\ite{$d_n$} {c.f.\ Definition \eqref{dN-Eqn}},
\ite{$e^{\widehat A}$} {$K$-theory Euler class corresponding to the
  $\widehat A$-genus},
\ite{$e^{\td}$} {$K$-theory Euler class corresponding to the
  $\td$-genus},
\ite{$\on{eq}(f,g)$} {equalizer of the maps $f$ and $g$},
\ite{$F_{(\underline g,\sigma)}$} {defined in Theorem \ref{Ls-Thm}},
\ite{$f_{(\underline g,\sigma)}$} {defined on page \pageref{fgs-page}},
\ite{$\mathcal G$} {one of the groups $\on{Spin(2n)}$ or $\on{U(n)}$},
\ite{$(\underline g,\sigma)$} {an element $\((g_1,\dots,g_n),\sigma\)$
  of the wreath product $G\wr\Sn$},
\ite{$K_{G,\Dev}$} {Devoto equivariant $K$-theory, c.f.\ Definition
  \ref{Dev-K-Def}},
\ite{$K_{G,\Dev,r}$} {Devoto equivariant $K$-theory for loops of
  length $r$, c.f.\ Definition \ref{Dev-K-Def}},
\ite{$K_{G}$} {$G$-equivariant $K$-theory},
\ite{${}_k\mL(M\mmod G)$} {orbifold loops of length $k$ in $M\mmod G$,
c.f. \ Definition \ref{kLMG-Def}},
\ite{${}_k\mV$} {the loop bundle $\mV$ rescaled, c.f.\ Definition
  \ref{kx-Def}}, 
\ite{${}_k(-)$} {see Defintion \ref{kx-ring-Def}},
\ite{$\Lambda_t$}{total exterior power}
\ite{$\Lambda(M\mmod G)$} {the inertia orbifold $\coprod_{[g]}M^g\mmod
  C_g$},
\ite{$M\mmod G$} {global quotient orbifold of an action of $G$ on $M$},
\ite{$M^g$} {the $g$-fixed points of $M$},
\ite{{$\on{MU}$}} {complex cobordism},
\ite{{$\on{MUP}$}} {periodic complex cobordism, c.f.\ Example
  \ref{MUP-def-Exa}}, 
\ite{{$\on{MSpinP}$}} {periodic spin cobordism, c.f.\ Example
  \ref{MUP-def-Exa}}, 
\ite{$N^g$} {normal bundle of $M^g$ in $M$},
\ite{$P_n = P^K_n$} {$n^{th}$ Atiyah power operation in $K$-theory},
\ite{$P^\str_n$} {$n^{th}$ stringy power operation on $K_\Dev$, c.f.\ Definition
  \ref{stringy-P-Def}},
\ite{$P^\top_n$} {$n^{th}$ Atiyah power operation on $K_G\ps q$, c.f.\
  Definition \ref{Atiyah-powerop-Def} or on $K_{\Dev}$, c.f.\ page
  \pageref{atiyah-powerop-in-Kdev-page}}, 
\ite{$P_t$} {total power operation},
\ite{$\pi_!$} {push-forward along the map $\pi$ (typically $\pi$ is
  the unique map to the one point space), c.f.\ \eqref{eq:pushforward-def}},
\ite{$q$} {the defining representation of the circle group $\mathbb
  R/\ZZ$},
\ite{$q^{\frac1l}$} {the representation of $\mathbb R/l\ZZ$ obtained
  by scaling $q$ with a factor $1/l$},
\ite{$q^{\frac jl}$} {the $j^{th}$ tensor power $\(q^{\frac1l}\)^{\tensor_\CC j}$},
\ite{$S_E(X)$} {(symmetric algebra) the graded ring
  $\bigoplus_{k\geq 0}E_{G\wr 
    \Sn}(X^n)t^n$, c.f.\ page \pageref{S_E-page}},
\ite{$\on{Sym}_t$} {total symmetric power},
\ite{$U^2$-bundle} {a vector bundle with compatible $\on{Spin}$- and
  complex structures, c.f.\ \eqref{U2-Eqn}}, 
\ite{$u^{\widehat A}$} {$K$-theory Thom class corresponding to the
  $\widehat A$-genus},
\ite{$u^{\td}$} {$K$-theory Thom class corresponding to the
  $\td$-genus},
\ite{$u^\str_{\widehat A}$} {stringy (equivariant) $\widehat A$-Thom
  class},
\ite{$u^\str_{\td}$} {stringy (equivariant) $\td$-Thom class},
\ite{$V_1$} {the summand of the bundle $V\at{M^g}$ fixed by $g$},
\ite{$(V_\CC)_\zeta$} {the summand of the complexified bundle
  ${V_\CC}\at{M^g}$ on which $g$ acts with eigenvalue $\zeta$},
\ite{$\underline X$} {an orbifold},
\ite{$X^V$} {Thom space of a vector bundle $V$ over a space $X$},
\ite{$\zeta_l$} {the primitive $l^{th}$ root of unity $e^{2\pi i/l}$},
\ite{$g\sim h$} {the group elements $g$ and $h$ are conjugate},
\ite{{$[g]$}} {conjugacy class of $g$},
\ite{$M\rtimes G$} {the translation groupoid of an action of $G$ on $M$},
\ite{$\boxtimes$} {external tensor product (different groups), c.f.\
  Definition \ref{products-def}}, 
\ite{$\tensor$} {depending on the context, external or internal tensor
  product over $G$}.
\end{list}

\section{Localization and Fourier expansion}\label{Fourier-Sec}
This section introduces some basic definitions. We start with the
spaces modelling loops of length $k$ in $M\mmod G$. 
Our definitions are a generalization of the familiar case $k=1$
(see for example \cite{FNLU}).
We need to include the case $k>1$ for technical reasons: it will
be necessary when we consider iterated symmetric powers.
For the time being the reader is invited to ignore this issue and to set
$k=1$. Throughout the
paper we will drop $k$ from the notation if $k=1$. 
All our paths are piecewise smooth, and path-spaces carry the
compact-open topology. 

\subsection{Orbifold loop spaces}
Let $G$ be a finite group acting smoothly from
the right on a manifold $M$. Let $g$ be an element of $G$, and let
$l=|g|$ be the order of
$g$. We write $C_g = C_G(g)$ for the centralizer of $g$ in $G$.
\begin{Def}
  For a natural number $k$, we define
    $${}_k\mP_gM := \{\map\gamma{[0,k]}M\mid\gamma(k)=\gamma(0)g\}$$
  and 
    $${}_k\mL_gM := \maps_{\ZZ/l\ZZ}(\RR/kl\ZZ,M)$$
  ($\ZZ/l\ZZ$-equivariant maps), where $1\in\ZZ/l\ZZ$ acts as $k$ on
  $\RR/kl$ and as $g$ on $M$. 
\end{Def}
The centralizer $C_g$ acts on both of these spaces via its action on
$M$, and there is a $C_g$-equivariant homeomorphism
\begin{eqnarray}\label{loops-paths-Eqn}
  {}_k\mathcal P_gM & \cong &
  {}_k\mL_gM\\
  \notag
  \gamma &\mapsto &\gamma*\gamma g*\dots*\gamma g^{l-1}.
\end{eqnarray}
  \begin{figure}[h]
  \label{fig:PL}
    \centering
\includegraphics[scale=.35]{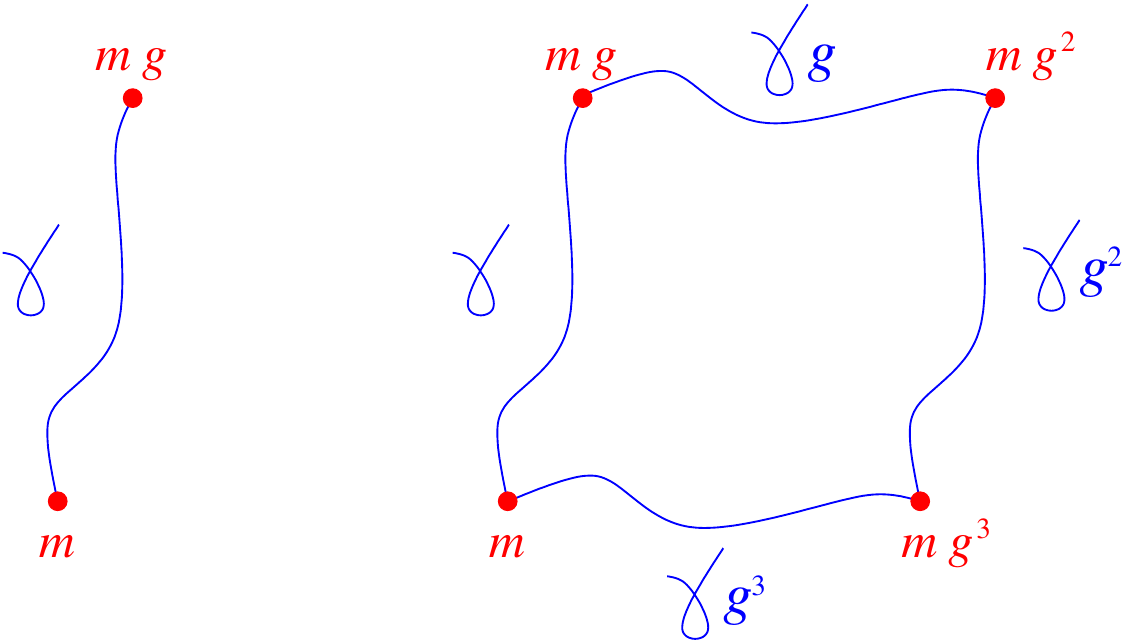}
    \caption{A path $\gamma\in\mP_gM$ and its image in $\mL_g M$. Here
      $l=4$.} 
\label{Path-Loop:Fig}
  \end{figure}

\begin{figure}[h]
    \centering
\includegraphics[scale=.35]{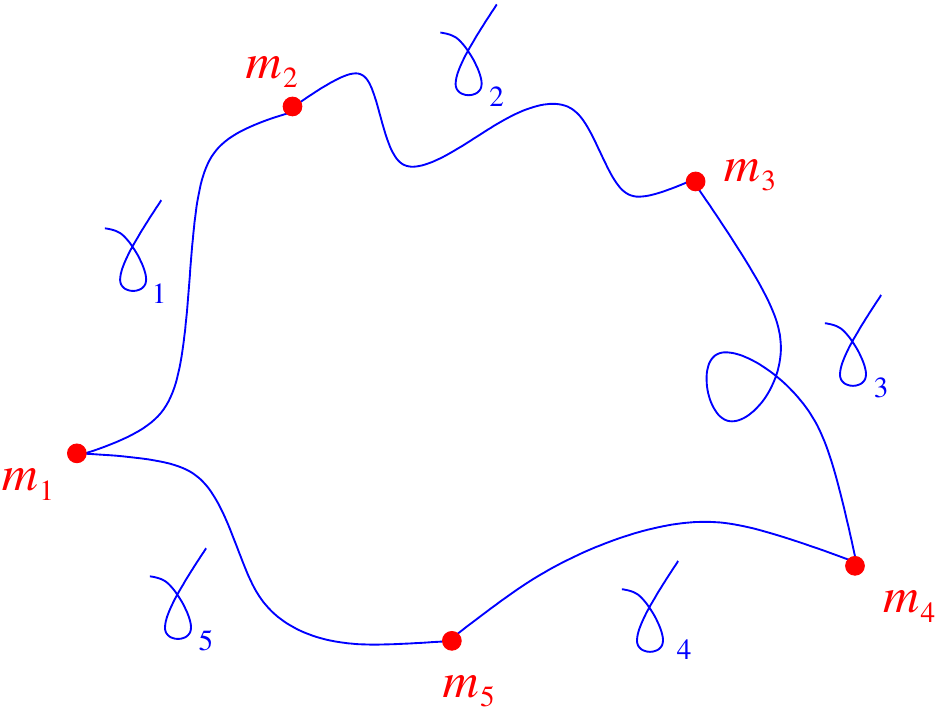}
    \caption{An element $\gamma = \gamma_1*\dots*\gamma_5\in{}_5\mL M$ for trivial $G$.}
\label{Path-of-length-5:Fig}
  \end{figure}

The group $\RR/kl\ZZ$ (and hence also its subgroup $\ZZ/kl\ZZ$) acts on
${}_k\mL_gM$ by rotation of the loops. Note that this action commutes
with that of $C_g$, and that the actions
of $k\in\RR/kl\ZZ$ and $g\in C_g$ agree. 
This motivates the following definition.
\begin{Def}\label{extended-centralizer-def}
  Let $k$ be a natural number. We write $C^k_g$ or $C_G^k(g)$ for the
  quotient of $C_g\times\ZZ/kl\ZZ$ by the normal subgroup generated by
  $(g\inv,k)$, and we write $\CRg$ for the quotient of
  $C_g\times\RR/kl\ZZ$ by the same normal subgroup. 
\end{Def}
It follows from the discussion above that we have actions of
$\CRg$ and $C^k_g$ on ${}_k\mL_gM$. Using the homeomorphism
\eqref{loops-paths-Eqn}, we 
obtain $\CRg$ and $C_g^k$-actions on ${}_k\mP_gM$.
Throughout the paper we will use the homeomorphism
\eqref{loops-paths-Eqn} to identify
${}_k\mL_gM$ with 
${}_k\mP_gM$ and, by abuse of notation, we will write ${}_k\mL_gM$ for
both.
\begin{Def}\label{kLMG-Def}
  Let ${}_k\mL(M\mmod G)$ be defined by
  $${}_k\mL(M\mmod G) := \coprod_{[g]}{}_k\mL_gM\mmod C^k_g,$$
  where the disjoint union is over the conjugacy classes of $G$.
\end{Def}
Think of ${}_k\mL(M\mmod G)$ as the groupoid of $k$ open strings
in $M\mmod G$ joining together to form one long closed string, where the
order of the $k$ strings does not matter (see Figure
\ref{Path-of-length-5:Fig}).  
\subsection{Representations of $\RR/l\ZZ$}
The goal of this section is to set the stage for our Fourier expansion
principle in Section \ref{Sec:loop}.
\begin{Def}
  We let 
  \begin{eqnarray*}
    q\negmedspace :\mathbb R/\ZZ &   \to   & U(1) \\
                              t  & \mapsto & e^{2\pi i t}
  \end{eqnarray*}
  denote the fundamental (complex) representation of the circle group 
  $\mathbb R/\ZZ$. More generally, we write
  \begin{eqnarray*}
    q^{\frac1l}\negmedspace :\mathbb R/\ZZ &   \to   & U(1) \\
                              t  & \mapsto & e^{2\pi i t/l}
  \end{eqnarray*}
  for the representation of the long circle $\mathbb R/l\ZZ$ obtained by
  scaling the exponent $t$ by $1/l$.
\end{Def}
Note that
$$1\in\ZZ/l\ZZ\sub\RR/l\ZZ$$ 
acts on $q^{\frac 1l}$ by multiplication with
a primitive 
$l^{th}$ root of unity, which we will denote $\zeta_l$.
Recall that all non-trivial irreducible real representations of $\Rl$ 
are of the form
$$q^\frac jl_\RR:=\res\CC\RR q^\frac jl,$$ 
where $\res\CC\RR \rho$ denotes
the underlying real 
vector space of $\rho$. Further, $q^\frac il\cong_\RR q^\frac jl$ if and
only if $i=\pm j$.
Let $\mG$ be one of the groups $U(n)$ or $\Spin(2n)$, and let
$$\rho\negmedspace :\Rl\to\mG\to SO(2n)$$ 
be a $\mG$-representation of
$\Rl$ with underlying vector space $V$. Then $V$ decomposes into a
direct sum 
\begin{equation}
  \label{decomp-Eqn}
  V\cong_\RR V_0 \oplus \bigoplus_{j\geq 1} V_jq^{\frac jl}.
\end{equation}
Here $V_jq^\frac jl$ denotes the summand with rotation number $j$,
and each of the summands is a complex/$\Spin$ representation.
The notation ``$V_jq^j$'' is motivated as follows:
for $j\geq 1$, we can endow the underlying vector space $V_j$ of
$V_jq^\frac jl$ with a complex structure, where 
multiplication with $i$ is defined as the action of $\frac l{4j}$.  
Then
$$V_jq^\frac jl = \res\CC\RR \(V_j\tensor_\CC \(q^\frac1l\)^{\tensor j}\).$$

{\bf Warning:} if the structure group $\mG$ is $U(n)$, then $V$ decomposes as a complex
representation
$$
  V\cong_\CC \bigoplus_{j\in\ZZ}W_j q^\frac jl.
$$
However, (\ref{decomp-Eqn}) endowed with the  
complex structure defined by the rotations is
\begin{equation}
  \label{cx-decomp-Eqn}
  V\cong_\RR W_0\oplus\bigoplus_{j\geq 1} \(\overline{W}_{-j}\oplus
  W_j\) q^\frac jl.
\end{equation}
\subsection{Loop bundles and Fourier decomposition}\label{Sec:loop}

We are now ready to define loop vector bundles over orbifold loop spaces.
The definition is somewhat technical; our main interest will be in the
motivating Example \ref{motivate-loop-bdl-Exa}.
\begin{Def}\label{S-bdl-Def}
  Let $\mG$ be one of the groups $U(n)$ or $\Spin(2n)$. Let $G$ be a
  finite group acting smoothly from the right on $M$, let $g\in G$
  have order $|g|=l$.
  An (orbifold) $\SS$-equivariant $\mL\mG$-vector bundle over ${}_k\mathcal
  L_gM\mmod C_g^k$ is a (typically infinite dimensional) vector bundle
  $\mV$ over ${}_k\mathcal L_gM$ 
  with structure group ${}_{kl}\LG$, 
  which is $C_g^{\RR/k}$-equivariant as ${}_{kl}\mL\mG$-bundle
  (i.e., the groups act by ${}_{kl}\mL\mG$-bundle automorphisms),  
  such that
  the action of $\RR/kl\ZZ$ on $\mV$ intertwines
  with its action on ${}_{kl}\LG$.
\end{Def}
\begin{Exa}\label{motivate-loop-bdl-Exa}
  Let $\VV\in\on{Vect}(\XX)$ be a $\mG$-vector
  bundle on the orbifold 
  $\XX$. 
  View $\VV$ as an orbifold over $\XX$. Then the loop orbifold $\LV$ of
  $\VV$ is an orbifold $\SS$-equivariant $\mL\mG$-vector bundle over
  $\LX$. 
\end{Exa}

\begin{Def}
Two ${}_l\mL\mG$-vector bundles $\mathcal V$ and $\mathcal W$ 
over a space $X$ are called
{\em densely isomorphic} if there are dense sub-bundles $\mathcal
V'\sub\mathcal V$ and $\mathcal W'\sub\mathcal W$ and an isomorphism
of vector bundles with structure group $\mathcal L\mG$ between $\mathcal
V'$ and $\mathcal W'$. If in addition $\mathcal V$
and $\mathcal W$ are $G$-equivariant bundles over a $G$-space $X$,
then $\mathcal
V$ is densely isomorphic to $\mathcal W$ as a loop bundle with this extra
structure, if the inclusions and the vector bundle isomorphisms
preserve this extra structure. We will denote dense isomorphisms by
``$\cong$'' and often drop the word ``{\em densely}''. 
\end{Def}

Let $\mathcal V$ be an orbifold $\SS$-equivariant
$\LG$ vector bundle over ${}_k\LgM\mmod\Ckg$.
The condition that the action of $\RR/kl\ZZ$ on $\mV$ intertwines with
that of ${}_{kl}\LG$ spells out to the following: the two actions combine to a
(left-) action of $\RR/kl\ZZ\ltimes_{kl}\LG$ on $\mV$, where 
$$
  (t,\gamma')\in \RR/kl\ZZ\ltimes{}_{kl}\LG
$$
maps $v\in\mV_\gamma$ to 
$$
  (t,\gamma')\cdot v = t\cdot (\gamma'\cdot v) = (t\cdot
  \gamma')(t\cdot v) \in \mathcal V_{t\cdot \gamma}.
$$
Restricting to constant loops, we get an action of $\CRg\times\mG$ on
$\mV\at{M^g}$, where $\RR/kl\ZZ\times\mG$ acts fiber preserving.
Therefore, the restriction $\mathcal V\at{M^g}$ decomposes into a direct sum of
$\RR/kl\ZZ$-representations
$$
  \mathcal V\at{M^g}\cong_\RR V_0\oplus \bigoplus_{j\geq 1}V_jq^\frac j{kl}  
$$
as in \eqref{decomp-Eqn} above,
where $V_jq^\frac j{kl}$ is a (finite dimensional) $\CRg$-equivariant
$\mG$-vector bundle over $M^g$.  
The actions of $g\in C_g$ and $k\in\RR/kl\ZZ$ on
$V_j$ agree, so that $g$ acts on $V_j$ as (complex) multiplication by
$\zeta_l^j$. 
\begin{Exa}
  Let $V$ be a finite dimensional real vector space, and recall (e.g.\
  from \cite[6]{Hirzebruch:Berger:Jung}) that there is a dense
  isomorphism of real $\Rl$-representations
  $$
    \maps(\Rl,V) \cong V \oplus\bigoplus_{j\geq1} V_\CC q^\frac jl,
  $$
  given by the Fourier expansion principle.  
  Here $V_\CC$ denotes the complexification of $V$.
  Let now $\gg$ act on $V$, where $|g|=l$. 
  Note that $1\in\Rl$ acts by multiplication with $\zeta_l^j$ on the
  $j^{th}$ summand of the right-hand side.
  A loop is $\ZZ/l\ZZ$-equivariant, if and only if its Fourier
  expansion $\sum v_jq^j$ satisfies 
  $$
    \(\forall j\geq  0\) \quad v_jg=\zeta_l^jv_j.
  $$
  Hence we have
  $$
    \maps_{\ZZ/l\ZZ}(\Rl,V) \cong
    V_1\oplus \bigoplus_{j\geq 1}(V_\CC)_{\zeta_l^j} q^{\frac jl}.
  $$
  Here $(V_\CC)_\zeta$ denotes the $\zeta$-eigenspace of the
  complexification of $V$, and similarly, $V_1$ denotes the
  vectors of $V$ fixed by $g$.
\end{Exa}
\begin{Exa}\label{LVG-Exa}
  Let $\mV=\mL(V\mmod G)$ be the orbifold loop bundle of a $G$-equivariant 
  $\mG$-vector bundle over $M$ as in Example
  \ref{motivate-loop-bdl-Exa}. Let $l$ be the order of $g\in G$, 
  and let $m$ be a point in $M^g$. We have
  $$
    \mL_gV\at m \cong \maps_{\ZZ/l\ZZ}\(\Rl,V\at m\)
  $$
  and hence
  $$
    \mL_gV\at{M^g}\cong V_1 \oplus \bigoplus_{j\geq
      1}(V_\CC)_{\zeta_l^j}q^{\frac jl},
  $$
  where $(V_\CC)_\zeta$ denotes the $\zeta$-eigenbundle of the action
  of $g$ on $V_\CC\at{M^g}$.
 (Compare also \cite[(5.1.1)]{FNLU}.)
\end{Exa}

\section{Thom classes}\label{Thom-Sec}
\subsection{Devoto's equivariant Tate $K$-theory}
Taking $g$-fixed points is a functor 
$$
  \map{(-)^g}{G-\on{spaces}}{C_g-\on{spaces}}
$$
which preserves cofibre sequences, wedges and weak equivalences. Therefore
\begin{equation}\label{K-Dev-Eqn}
  X\mapsto \bigoplus_{[g]}K_{C_g}(X^g)\ps{q^\frac1{|g|}}    
\end{equation}
satisfies the axioms of a $G$-equivariant cohomology theory in
\cite{May:Alaska}. 
If $g$ is conjugate to $h$ in $G$ (denoted $g\sim h$) there is a
canonical isomorphism 
$$K_{C_g}(X^g)\cong K_{C_h}(X^h).$$
Thus (\ref{K-Dev-Eqn}) is, up to canonical isomorphism, independent of
choices of the representatives.
Example \ref{LVG-Exa} motivates the following definition.
\begin{Def}\label{Dev-K-Def}
  Consider the ring
  $$
    \bigoplus_{[g]}K_{C_g^r}(X^g)\ps{q^\frac1{r|g|}}.
  $$
  Note that in the $[g]$-summand, the coefficient $a^{[g]}_j$ of $q^{\frac
    j{r|g|}}$ is a virtual $C_g^r$-equivariant vector bundle. 
  The ring $K_{\Dev, G,r}(X)$ is defined as the subring of those power
  series where $$1\in \ZZ/r|g|\ZZ\sub C_g^r$$ acts by complex multiplication 
  with $\zeta_{r|g|}^j$ on the coefficient $a_j^{[g]}$.
  We will refer to this ring as the  
  {\em Devoto equivariant Tate $K$-theory ring (for strings of
  length $r$) of $X$}.
  As usual, if $r=1$, we will drop it from the notation. 
\end{Def}
\begin{Rem}
For $r=1$, the condition of the definition should be compared to Condition (a) of
the generalized Moonshine axioms (c.f.\ \cite{Norton}): the M\"obius
transformation 
$\tau\mapsto \tau+1$ (sending $q^\frac j{|g|}$ to $q^\frac
j{|g|}\zeta_{|g|}^j$) should correspond to the transition from $(g,h)$
to $(g,gh)$, i.e., to the action of $g$ on the summand corresponding to
$[g]$.  
\end{Rem}
%
%
\begin{Def}\label{products-def}
  Let $\map\alpha HG$ be a map of finite groups, and let
  $x\in K_{\Dev, G}(X)$. Then 
  $$
    \(\res{}\alpha x\)_{[h]} := x_{[\alpha(h)]}.
  $$

  Let $x\in K_{\Dev, G}(X)$ and $y\in K_{\Dev, H}(Y)$. Then the
  external tensor product of $x$ and $y$ is the element $$x\boxtimes
  y\in K_{\Dev, G\times H}(X\times Y)$$ defined by 
  $$
    \(x\boxtimes y\)_{[g,h]} := x_{[g]}\boxtimes y_{[h]}\in
    K_{C_G(g)\times C_H(h)}\(X^g\times Y^h\)\ps{q^\frac{1}{N}},
  $$ 
  where $N$ is the smallest common multiple of $|g|$ and $|h|$,
  $\boxtimes$ is the external tensor product of vector bundles, and
  $q^r\boxtimes q^s := q^{r+s}$. 
  If $G=H$, the $G$-equivariant external tensor product is defined as
  $$x\boxtimes y := \res{}\delta (x\boxtimes y),$$
  where $\delta$ denotes the diagonal inclusion of $G$ in $G\times
  G$. 
  If in addition $X=Y$, then the internal tensor product of $x$ and
  $y$ is defined to by the pullback along the diagonal $d$ of $X\times
  X$ of $x\boxtimes y$:
  $$
    x\tensor y := d^*(x\boxtimes y).
  $$ 
\end{Def}
%
%
%
%
%
%
%
\subsection{Stringy orbifold Thom classes}
This section serves as a morivation for Definition
\ref{stringy-A-hat-Def}. Much of the discussion here is based on the
circle of ideas described in the appendix. 
Let $G$ be a finite group, let $g$ be an element of $G$, and let $l$
be the order of $g$. Let $V$ be a $G$-equivariant $\Spin(2d)$-vector
bundle over a $G$-space $X$.
By Example \ref{LVG-Exa}, we have
\begin{equation}
  \label{str-decomp-Eqn}
  \mL_gV\at{X^g}\cong V_1\oplus \bigoplus_{j\geq 1} (V_\CC)_{\zeta_l^j} q^\frac jl.
\end{equation}

Assume that $X^g$ is connected.
The bundle $V_1$ is the bundle of $g$-fixed points of $V\at{X^g}$, and
the $g$-fixed points of the Thom space of $V$ are, as a $C_g$-space,
homeomorphic to the Thom space $\(X^g\)^{V_1}$ of $V_1$.

Assume first that $g$ acts trivially on $V\at{X^g}$, then $V_1=V\at{X^g}$, and
(\ref{str-decomp-Eqn}) contains only the terms with integral powers
of $q$. This case was treated in \cite{Haynes:loops}, we shortly recall the
discussion there:
Since $V$ is oriented, $V_\CC$ is an $\on{SU}$-vector
bundle. Hence it has a $\Spin(4d)$ structure as well as a complex
structure with trivial determinant bundle
$$\det V_\CC = 1\in K_{C_g}(X^g),$$ and we have have an equality of
the $\widehat A$-Euler class and the Todd genus Euler class 
$$
  e_{C_g}^{\hat A}(V_\CC\at{X^g}) =   e_{C_g}^{\td}(V_\CC\at{X^g}) 
$$
(see appendix).
Formula (\ref{Todd-A-hat-Eqn}) applied to the $\Rl$-action
suggests to define the 
$C_g\times\Rl$-equivarinat Thom class of $\mL_gV\at{X^g}$ as

$$
  u^{\hat A}_{C_g}(V) 
    \cdot\prod_{j\geq
    1}q^{-d j}\cdot\Lambda_{-q^j}V_\CC.
$$
Here $  u^{\hat A}$ denotes the $\widehat A$-Thom class, and
the exterior powers are exterior powers of
$C_g$-representations.  
It is customary (c.f.\ \cite{Segal:Bourbaki}) to use a $\zeta$-function
renormalization to 
replace the divergent sum $\sum_{j>0}j$ in the exponent of $q$ with
$\zeta(-1) = -1/12$, and to define
$$
  u^{\str}_g(V) :=
  q^{\frac d{12}}\cdot u^{\hat
    A}_{C_g}(V)\cdot\prod_{j=1}^\infty\Lambda_{-q^j}V_\CC. 
$$
Next, assume that $g$ acts fixed point free on $V\at{X^g}$ and that
$-1$ is not an eigenvalue of $g$.
For $j=1,\dots,\lfloor\frac l2\rfloor$, let $\theta_j = 2\pi\frac jl$,
and let $V_j$ be the subspace of $V$ on which $g$ acts by rotations
with angle $\theta_j$. We equip $V$ with the unique complex structure
such that $g$ acts on $V_j$ by multiplication with $e^{i\theta_j}$.
Then $V_j\tensor\CC$ is of the form
$$V_j\tensor\CC\cong V_j\oplus \overline{V_j}.$$
Setting $V_{-j}:= \overline{V_j}$, we have
$$
  \mL V\at{X^g}\cong\bigoplus_{j=1}^{\lfloor\frac l2\rfloor}V_jq^\frac jl\oplus
\bigoplus_{k=1}^{\infty}\bigoplus_{j=-\lfloor\frac
  l2\rfloor}^{\lfloor\frac l2\rfloor} V_j q^{k+\frac jl}.
$$
For fixed $k$, 
\begin{eqnarray*}
\det\(\bigoplus_{j=-\lfloor\frac l2\rfloor}^{\lfloor\frac l2\rfloor}
V_j q^{k+\frac jl}\)&\cong&
\bigotimes_{j=-\lfloor\frac l2\rfloor}^{\lfloor\frac l2\rfloor}
q^{(k+\frac jl)d_j}\cdot \det( V_j)\\
&\cong&  
q^{2dk}\cdot \det (V_\CC).
\end{eqnarray*}
Here $d_j=\dim_\CC V_j$, and
$$
  \sum_{j=-\lfloor\frac l2\rfloor}^{\lfloor\frac l2\rfloor}
  d_j(k+\frac jl) = 2k(\sum_{j=1}^{\lfloor\frac l2\rfloor}d_j) =
  2dk,
$$
while $\det V_\CC$ is trivial.
In a similar manner, the determinant line bundle of the first
$\lfloor\frac l2\rfloor$ summands becomes
$$
  q^{\frac{\on{age}(g)}2}\cdot\det(V)
$$
with
\label{age-page}
$$
  \on{age}(g) = 2\cdot\sum_{j=1}^{\lfloor\frac l2\rfloor}\frac jld_j.
$$ 
Note that $\det (V)$ is not trivial, but has a square root.

As above, Formula (\ref{Todd-A-hat-Eqn})
suggests to define the $C_g\times\Rl$-equivariant Euler/Thom class of
$\mL_gV\at{X^g}$ as
$$
  q^{\frac d{12}-\frac{\on{age(g)}}{2}}\(\sqrt{\det V}\)\inv\cdot
  \prod_{j=1}^\infty\Lambda_{-\qjl}\(\(V_\CC\)_\zlj\).
$$

If we want to allow the eigenvalue $-1$, we need to assume that 
the eigenspace $V_{\frac l2}$ of the action of $-1$ on $V$ has a
complex structure. 
Then we can make sense of the expression $(\sqrt{\det V})$, and the
above discussion goes through.

This is for example the case if 
$V$ is a $G$-equivariant $U^2$-bundle (i.e., if $V$ has a complex and
a spin-structure that are compatible, see \eqref{U2-Eqn}).
Note that in that case we have two complex
structures on the other $V_j$  -- the given complex structure and  
the auxiliary one we defined in formula (\ref{decomp-Eqn}) above.
We choose to work with the latter. This is an issue similar to 
the one in the warning after Formula (\ref{decomp-Eqn}). 
%
%
%
%
%
Write $V_1$ for the subspace of $V\at{X^g}$
that is fixed by the action of $g$, write $N^g$ for its 
orthogonal complement, and note
that 
$$
  u^{\hat A}_{C_g}(V_1)(\sqrt{\det N^g})\inv =
  u^\td_{C_g}(V_1)(\sqrt{\det V_1})\inv(\sqrt{\det N^g})\inv =
  u^{\td}_{C_g}(V_1)(\sqrt{\det V})\inv.
$$
All of the above discussion taken together motivates the following definition:
\begin{Def}\label{stringy-A-hat-Def}
  Let $V$ be a $G$-equivariant $U^2(d)$-vector bundle over a
  compact $G$-space $X$. 
  Let $g$ be an element of $G$ with order $l$.
  Assume that $X^g$ is connected.
  We define
  $$
    u_g^\str(V) := q^{\frac{d}{12}-\frac{\on{age}(g)}{2}} \cdot
       u^{\td}_{C_g}(V_1)\cdot (\sqrt{\det V})\inv
    \cdot\prod_{j=1}^\infty\Lambda_{-\qjl}\(\(V_\CC\)_\zlj\)
  $$
  and set
  $$
    u^\str(V) := \sum_{[g]}u^\str_g(V).
  $$ 
\end{Def}
This is an element of $K_{\Dev,G}\(X^V\)[q^{\pm\frac 1{12}}]$.
If $X^g$ is not connected, one sums over its connected components.
\begin{Prop}\label{ustrings-are-thom-classes-Prop}
The $u^\str_G(V)$ 
  satisfy the axioms of Definition \ref{Thom-Def}.
\end{Prop}
\begin{Pf}{}
  The naturality axiom and the change of groups axiom are
  straight-forward. Multiplicativity follows from that of $u^\td$ 
  and from the multiplicative properties of $\sqrt{\det(-)}$ and of
  the total exterior power $\Lambda_t$. Further, $\Lambda_{-t}$ is
  invertible with inverse $\on{Sym}_t(V)$, which implies the
  periodicity axiom.   
\end{Pf}

\subsubsection{Comparison to the work of Kitchloo and Morava}
Let $G$ be the trivial group, and let $V$ be a complex vector bundle
of dimension $d$. 
Using (\ref{cx-decomp-Eqn})
and the fact that
\begin{eqnarray*}
  \Lambda_{-1}\(Vq^{-k}\oplus Vq^k\) &\cong& \det\(-Vq^{-k}\)\cdot  
  \Lambda_{-1}\(V_\CC q^k\)\\
  &\cong& q^{dk}(\det V)\inv \cdot\Lambda_{-1}\(V_\CC q^k\),
\end{eqnarray*}
we see that if $\det(V)$ is trivial, a similar discussion as
above suggests the definition
$$
  u_\td^\str(V) := q^{-\frac1{12}}u_\td(V)\cdot\prod_{k\geq 1}\Lambda_{-q^k}V_\CC.
$$
The periodicity element
$\alpha=u^\str_\td(\CC)$ (c.f.\ (\ref{alpha})) is 
$$
  \alpha = \(q^{-\frac1{12}}\cdot\prod_{j\geq 1}(1-q^j)^2\)\cdot\beta,
$$
where $\beta$ is the Bott-element.
Hence the remormalized Thom class in degree $2d$
\begin{eqnarray*}
  \tau^\str_\td(V) & := & u^\str_\td(V)\cdot\alpha^{-d} \\
  &=& \tau^{\td}(V)\cdot\prod_{j\geq 1}\frac{\Lambda_{-q^j} V_\CC}{(1-q^j)^{2d}}  
\end{eqnarray*}
is the Thom class defined in \cite{Ando:Morava}, \cite[4.1]{Kitchloo:Morava}. 
Here $\tau^{\td}(V) = u_{\td}(V)\cdot\beta^{-d}$.
\subsection{The equivariant Witten genus}
Along with the existence of equivariant Thom classes comes an
equivariant Witten genus. We are now going to calculate this
equivariant Witten genus in terms of characteristic classes.
\begin{Cor}[{of Proposition \ref{ustrings-are-thom-classes-Prop}}]\label{genus-Cor}
  There is a unique map of $G$-spectra
  $$
    \map{\phi_G}{\MU^2_G}{K_{\Dev,G}}
  $$
  taking Thom classes to Thom classes.
  The induced genus 
  $$\AA^\str_G(M):=\phi_G(M)$$ of a 
  $G$-equivariant $U^2$ manifold 
  $M$ of dimension $2n$ is given by the formula 
  $$
    \(\AA^\str_G(M)\)_{[g]} = 
    q^{-\frac{d}{12}+\frac{\ageg}2}\cdot
   \td_{C_g}\(M^g;\sqrt{(TM)} \cdot 
    \prod_{j\geq 1} 
    \on{Sym}_{q^\frac jl} \((TM^\CC)_{\zeta_l^j}\)\),
  $$
where $l$ is the order of $g$, and we sum over the connected
components of $M^g$.  
\end{Cor}

\begin{Pf}{}
  The existence and uniqueness of $\phi_G$ follow from Proposition
  \ref{Okonek-Prop}. Further, assuming $M^g$ connected, we have
  \begin{eqnarray*}
    \(\AA^\str_G(M)\)_{[g]} &\cong &\(\pi_!^\str(1)\)_{[g]}\\
    &\cong & \pi_{C_g!}^\td\(\(q^{\frac{d}{12}-\frac\ageg2}(\sqrt{TM})\inv\prod_{j\geq 1}\Lambda_{-q^{\frac
        jl}}(TM^\CC)_{\zeta_l^j}\)\inv\)\\ 
    &\cong &
    q^{-\frac{d}{12}+\frac\ageg2}\cdot\td_{C_g}\(M^g,\sqrt{TM}\prod_{j\geq
      1}\on{Sym}_{q^{\frac 
        jl}}(TM^\CC)_{\zeta_l^j}\),
  \end{eqnarray*}
  where the first and last equation are (\ref{push-forward-Eqn}), and the
  second equation follows from Definition \ref{stringy-A-hat-Def}.
\end{Pf}

It is traditional to express equivariant genera in the following more
explicit and somewhat more complicated form.
\begin{Cor}\label{Astring-gh-Cor} 
  Let $h\in C_g$, and $k=|h|$. Assume $M^g$ to be connected. Then we have
  $$
    \AA^\str_G(M)_{[g]}(h) =
    \pm q^{-\frac{d}{12}+\frac\ageg2} e^{i\pi\medspace\ageh}
     \int_{M^{g,h}} e^H\(M^{g,h}\) \cdot 
    \(\prod_{\tilde x_i}e^{\frac{\tilde x_i}{2}}\)$$ $$\cdot
    \left[ 
    \prod_{s,x_i}\(1-\zeta_k^s e^{x_i}\)
  \prod_{j\geq 1}\prod_{s,y_i}
  \(1-\zeta_k^se^{y_i}q^\frac jl\)
    \right]\inv,
  $$
  where $e^H$ denotes
  the Euler class in ordinary cohomology,
  $s$ runs from $0$ to $k-1$, 
  the $\tilde x_i$ run over the Chern roots of $TM$, 
  the $x_i$ are the Chern roots of
  $TM_{g=1,h=\zeta_k^s}$, and
  the $y_i$ are the Chern roots of the simultaneous eigenspaces 
  $(TM_\CC)_{g=\zeta_l^j,h=\zeta_k^s}$.
  (All bundles are over the double fixed points $M^{g,h}$.)  
\end{Cor}
\begin{Pf}{}
By (\ref{Lefschetz-Eqn}), we have
\begin{equation}\label{computation-of-Astring-Eqn}
    \td_{C_g}\(M^g;a\)(h) = 
    \int_{M^{g,h}} {e^H\(M^{g,h}\)\cdot\on{ch}(a(h))\cdot
    \prod_{s,x_i} \frac1
    {1-\zeta_{k}^{s} e^{{x_i}}}
    },
\end{equation}
where $s$ runs from $0$ to $k-1$, and the $x_i$ are the Chern roots of
$TM_{g=1,h=\zeta_k^s}$. 
We have
$$
  \on{ch}(\sqrt{TM}(h)) = \pm \prod_{s,\tilde x_i}\zeta_{2k}^s
  e^{\frac{\tilde x_i}{2}}=
  \pm e^{{i\pi\medspace\ageh}}\cdot\prod_{\tilde x_i}e^{\frac{\tilde x_i}{2}},  
$$
where the $\tilde x_i$ are the Chern roots of $(TM)_{h=\zeta_k^j}$.
The Chern character of the remaining factor, evaluated at $h$, computes to
$$
  \prod_{j\geq 1}\prod_{s=0}^{k-1}\prod_{y_i}
  \(1-\zeta_k^se^{y_i}q^\frac jl\),
$$
where the $y_i$ are the Chern roots of the simultaneous eigenspaces
$
  (TM_\CC)_{g=\zeta_l^j,h=\zeta_k^s}.
$
\end{Pf}

Corollary \ref{Astring-gh-Cor} should be compared to the calculations
on the last page of \cite{FNLU}. Basically, our discussion starts, once
the authors of \cite{FNLU} replace their variable $\hat q$ with its monodromy $2\pi i$.
For our purposes, the double loop space picture is not necessary. It
will, however, become very important when one wants to study
modularity properties, i.e., the ``moonshine-like'' behavior discussed
in \cite{Devoto}: 
We already observed that replacing the pair of commuting elements
$(g,h)$ with $(g,gh)$ corresponds to replacing 
$\tau$ with $\tau+1$, where $q=e^{2\pi i\tau}$.
This is half of the first Moonshine condition, the other half demands
that switching to $(h\inv,g)$ should correspond to $\tau\mapsto
-\frac1\tau$. In the picture described in the beautiful paper
\cite{FNLU}, the elements $g$ and $h$ describe
the monodromy of a $G$-principal bundle $P$
over a torus, where $P$ maps equivariantly to $M$.
As described in \cite{Freed:Quinn}, \cite{Ganter:twistedhecke}, the transformation
behaviour of generalized Moonshine corresponds to a
different choice of circles, generating the same torus. 
Our approach sees the first circle, but not the second one. This is
the reason why we can explain the first transformation geometrically,
but not the second one.
%
%
\begin{Def}\label{orbifold-Witten-Def}
  Let $M$ be a $G$-equivariant $U^2(d)$-manifold. Then 
  the orbifold Witten genus of $M\mmod G$ is defined by
  $$
    \AA^\str_\orb(M\mmod G) :=\GG G\sum_{gh=hg}\(\AA^\str_G(M)_{[g]}\)(h),
  $$
  where the sum is over all pairs of commuting elements of $G$.
\end{Def}
One of the advantages of working with our setup is that the proofs of
integrality and Morita invariance become fairly straightforward.
\subsection{Integrality}\label{Integrality-Sec}
\begin{Prop} 
  The orbifold genus $\AA_\orb^\str$ takes values in
  $\ZZ\ps q[q^{\pm\frac 1{12}}]$. 
\end{Prop}
\begin{Pf}{} Let $g\in G$ have order $l$. Then 
  $\AA^\str_G(M)_{[g]}$ is an element of
  $R(C_g)\ps{q^\frac1l}[q^{\pm\frac1{12}}]$, and 
  \begin{eqnarray*}
    \AA^\str_\orb(M\mmod G) & = & 
    \sum_{[g]}\GG{C_g}\sum_{h\in C_g}\AA^\str_G(M)_{[g]}(h) \\
    & = & \sum_{[g]}\langle\AA^\str_G(M)_{[g]},1\rangle_{C_g}.
  \end{eqnarray*}
  Here $\langle-,-\rangle_{C_g}$ denotes the inner product in $R(C_g)$,
  so that
  \begin{eqnarray*}
    \left\langle\sum V_j q^\frac jl,1\right\rangle_{C_g} &:=& 
    \sum \langle V_j,1\rangle_{C_g}  q^\frac jl\\    
    & = & \sum\dim_\CC\(V_j^{C_g}\)q^\frac jl
  \end{eqnarray*}
  computes the dimensions of the maximal trivial summands of the $V_j$.
  Fix $g$ for the remainder of the proof, and assume, for simplicity,
  that $M^g$ is connected.
  The element $g\in C_g$ acts by multiplication with
  $\zeta_l^k$ on the coefficient of $q^\frac kl$ in 
  $$\prod_{j\geq 1}\Lambda_{-q^{\frac jl}}(TM^\CC)_{\zeta_l^j}.$$
  Further, $g$ acts by multiplication with $\pm e^{-i\pi\medskip\ageg}$ on
  $$
    q^{-\frac\ageg2} u_\td(TM^g)\cdot(\sqrt{\det TM})\inv.
  $$
  To determine the sign, note that $g$ acts trivially on $u_\td(M^g)$ and
  by multiplication with
  $e^{2i\pi\medskip\ageg}$ on $\det TM$, and that the sign 
  depends on the choice of an identification of $\Spin(2)$ with $U(1)$.
  For our calculation of $\AA^\str(M)_{[g]}$, we
  fix this identification in such a way that the sign becomes positive.
  Then only the coefficients of integral powers of $q$ in
  $$
    \pi^{M^g}_{\td,!}\(\(q^{-\frac\ageg2}\cdot \sqrt{\det TM}\inv\cdot
    \prod_{j\geq 1}\Lambda_{-q^{\frac jl}}(TM^\CC)_{\zeta_l^j}\)\inv
    \)
  $$
  contribute
  trivial summands.
  If we had chosen to work with the other choice of identification, we
  would have to multiply the above result for
  $\AA^\str(M)_{[g]}$ with $-1$, which would not destroy the integrality. 
\end{Pf}
\subsection{The Atiyah-Segal map for orbifolds and Morita invariance
  of the orbifold Witten genus}
A priori, it is not clear that Definition \ref{orbifold-Witten-Def} is
independent of the presentation of the orbifold as a global quotient
$M\mmod G$. For instance, why should an isomorphism of orbifolds
$M\mmod G\cong N\mmod H$ imply $$\AA^\str_\orb(M\mmod G) =
\AA^\str_\orb(N\mmod H)?$$  
This section proves this fact for complex orbifolds with a holomorphic
root of the line bundle $\det(T\underline X)$. This condition is, for
instance, satisfied by Calabi-Yau orbifolds.
The main tool of the proof is the Atiyah-Segal character map for
orbifolds. This map was defined in \cite{Adem:Ruan} and later in
\cite{Moerdijk}, but we will use a definition that is somewhat more
elementary and closely follows the original approach by Atiyah and
Segal \cite{Segal:handwritten}.

For the basic definitions of orbifolds, we refer the reader to the
introductory paper \cite{Moerdijk}. 
Let ${\bf G}$  be an orbifold groupoid, and recall that the inertia
groupoid $\Lambda({\bf G})$ as objects and morphisms
$$
  \Lambda(\bfG)_0 =\eq(s,t) \quad\text{ and }\quad \Lambda(\bfG)_1 =
  \eq(s,t) \times_{G_0}G_1.
$$
Here $\eq(s,t)$ stands for the equalizer of the source and target of $\bfG$),
the map $G_1\to G_0$ defining the fibred product is that target map
$t$, and $(g,h)\in\Lambda(\bfG)_1$ is an arrow from $g$ to $h\inv gh$.

Consider the map
\begin{eqnarray*}
  {A}\negmedspace : {\eq(s,t)}&\longrightarrow&
  \eq(s,t) \times_{G_0}G_1   \\
  g&\longmapsto& (g,g).
\end{eqnarray*}
This map $A$ is simultaneously a section of $s$ and $t$, and
one checks that $A$ is a natural transformation from the identity map of
$\Lambda(\bfG)$ to itself. 

Let $V$ be a vector bundle over the
groupoid $\Lambda(\bfG)$. By definition, $V$ is a vector bundle over
$\Lambda(\bfG)_0$ together with an isomorphism 
$$
  \map\mu{t^*V}{s^*V}
$$
over $\Lambda(\bfG)_1$ such that a few diagrams commute.
The section $A$ pulls back $\mu$ to a (fiber preserving) automorphism of
$$V=A^*t^*V = A^*s^* V$$
over $\Lambda(\bfG)_0$,
which turns out to be an automorphism of $V$ as orbifold vector
bundle. All of the above is preserved by equivalences of
orbifold groupoids $\bfG\simeq{\bf H}$.
\begin{Exa}\label{translation-Exa}
  Let $\bfG=M\rtimes G$ be the translation groupoid of the action of a finite
  group $G$ on a manifold $M$, and let $V$ be a $G$-equivariant vector
  bundle on $M$. 
  Then $\Lambda(\bfG)$ is equivalent to
  the disjoint 
  union over the conjugacy classes of $G$ of the translation groupoids
  of the actions of the centralizers on the fixed point loci:
  $$
    \Lambda(\bfG)\simeq \coprod_{[g]}M^g\rtimes C_g.
  $$
  Under this equivalence, the map $A$ corresponds to
  the maps
  \begin{eqnarray*}
    M^g &\to& M^g\times C_g\\
    m   &\mapsto & (m,g),
  \end{eqnarray*}
  and $A^*(\mu)\at{M^g}$ becomes multiplication with $g$ on
  $V\at{M^g}$. 
\end{Exa}
\begin{Def}
  Let $\underline X$ be an orbifold, and let $\Lambda(\underline X)$ be its inertia
  orbifold. We define
  \begin{eqnarray*}
    \chi'\negmedspace: K_\orb(\Lambda(\underline X)) &\to &
    K_\orb(\Lambda(\underline X))\tensor \CC\\
    V&\mapsto&\sum_\zeta V_\zeta\tensor \zeta,
  \end{eqnarray*}
  where the sum is over all eigenvalues of the action of $A^*(\mu)$ on
  $V$. Write $i$ for the canonical map
  \begin{eqnarray*}
    i:\Lambda(\underline X) &\longrightarrow &\underline X\\
           s=t : \eq(s,t) &\longrightarrow & G_0\\
           \on{pr_2} : \eq(s,t)\times_{G_0}G_1 &\longrightarrow & G_1.
  \end{eqnarray*}
  We define the {\em Atiyah-Segal character map} as the composite of
  $i^*$ with $\chi'$.
\end{Def}
\begin{Exa}
  In the situation of Example \ref{translation-Exa}, the map $i$
  becomes the inclusion of the fixed point locus in $M$, and 
  $$
    \chi\negmedspace : K_G(M)\longrightarrow
    \bigoplus_{[g]}K_{C_g}(M)\tensor \CC
  $$
  is the classical Atiyah-Segal character map.
\end{Exa}
Let now $\XX$ be a compact complex orbifold and $T\XX$ its tangent
bundle. 
It follows from Example \ref{translation-Exa} and the analogous
statement for the equivariant case that locally (and hence globally)
we have
$
  (T\XX)_1 \cong T\Lambda\XX.
$
Let $L$ be the smallest common multiple of the orders of all the
elements of the stabilizers of $\XX$, and recall that the connected
components of an orbifold are well-defined (they correspond to the connected
components of its quotient space). Further, it is still true for
orbifolds that the
dimension of a vector bundle is constant over each connected component.
Hence, over each connected component of $\XX$,  
$$
 q^{-\frac{d}{12}+\frac{\on{age}(A^*(\mu))}2}
 \prod_{j\geq 1}\on{Sym}_{q^\frac
  jL}\(T\XX^\CC\)_{\zeta_L^j} 
$$
is a well defined element of $\Korb{\Lambda\XX}\ps{q^\frac1L}[q^{\pm\frac1{12}}]$.
Here $T\XX^\CC$ is restricted to $(\Lambda\XX)$ along $i$.

Let $\mA$ be an abelian sheaf over the orbifold $\XX$. Recall that the
orbifold cohomology of $\XX$ with coefficients in $\mA$, denoted
$$H^*\(\XX,\mA\),$$
is defined as the cohomology of the complex
$\Gamma_{\on{inv}}(\mI^\bullet)$, where
$$
  0\to\mA\to\mI^\bullet
$$
is an injective resolution of $\mA$ in the category of abelian sheaves
over $\XX$ and $\Gamma_{\on{inv}}$ denotes the invariant sections.
\begin{Exa}\label{OV-Exa}
  Let $G$ be a finite group, $M$ a complex $G$-manifold, 
  $V$ a holomorphic $G$-vector bundle over
  $M$. We will write $\mathscr O(V)$ for the sheaf of germs of
  holomorphic sections of $V$.
  Let $\XX=M\mmod G$.
  Since
  taking $G$-invariants is an exact functor on complex vector spaces
  (denoted $(-)^G$), we have
  $$
    H^*(\XX,\mathscr O(V)) \cong  H^*(M,\mathscr O(V))^G.
  $$
\end{Exa}
\begin{Def}
  Let $\XX$ be a complex orbifold, and let $V$ be a holomorphic vector bundle on
  $\XX$. We define the topological Todd genus (also known as
  topological complex
  Euler characteristic) of $\XX$ with coefficients in $V$ by
  $$
    \td_\top(\XX,V) :=\sdim H^*(\XX,\mathscr O(V)) =
    \sum_p(-1)^pH^p(\XX,\mathscr O(V)).
  $$
  If $W$ is a holomorphic vector bundle over $\Lambda\XX$, then the
  orbifold Todd genus of $\XX$ with coefficients in $W$ is definied by
  $$
    \td_\orb(\XX,W) :=\sdim H^*(\Lambda\XX,\mathscr O(W)).
  $$  
\end{Def}
\begin{Exa}
  In the situation of Example \ref{OV-Exa}, one has
  \cite[pp.543-545]{Atiyah:Segal}
  \begin{eqnarray*}
    \td_\top(M\mmod G;V) & = & \sdim H^*(M;\mathscr O(V))^G\\
    &=& \GG G\sum_{g\in G}\sum_p(-1)^p\tr{H^p\(M,\mathscr O(V)\)}\\
    &=& \GG G\sum_{g\in G}\td\(M^g,\frac{V\at{M^g}}{\lambda_{-1}((N^g)^*)}(g)\).
  \end{eqnarray*}
\end{Exa}
\begin{Def} Let $\XX$ be a complex orbifold.
  We set
  $$
    \td_\top\(\XX;\sum_jV_j{q^\frac jl}\) :=  \sum_j\td_\top(\XX;V_j){q^\frac jl}.
  $$
  (here the $V_j$ are assumed holomorphic). 
  Assume that $\Lambda\XX$ is
  connected and that there exists a holomorphic line bundle
  $\sqrt{\det(T\XX)}$ on $\XX$ whose square is the determinant line
  bundle of the holomorphic tangent bundle $T\XX$. Then we set
%
  $$
    \AA^\str_{\orb}(\XX) := 
    q^{-\frac{d}{12}+\frac{\on{age}(A^*\mu)}2}\cdot \td_\orb\(\XX; 
    \sqrt{\det{T\XX}}\cdot
    \prod_{j\geq1}\on{Sym}_{q^\frac 
     jL}\(T\XX^\CC\)_{\zeta_L^j}\).     
  $$
  If $\Lambda\XX$ is not connected, one sums over its connected
  components. 
\end{Def}
In the case of a global quotient orbifold by a finite group, this
definition specializes to Definition \ref{stringy-A-hat-Def}.
\begin{Exa}
  An important class of orbifolds for which the conditions of the
  definition are satisfied are Calabi-Yau orbifolds: for these, the
  bundle $\det(T\XX)$ is trivial. 
\end{Exa}

The author recently learned of the work of
Dong, Liu and Ma, who defined Morita-invariant versions of
elliptic genera in \cite{Dong:Liu:Ma}. It seems likely that their
definitions are closely related to ours. A comparison to their work
would be very interesting, because they are working with the index
theorem on orbifolds. 
The authors of \cite{FNLU} also 
started from an index theoretic discussion but do not formulate their
results in a way that makes the Morita invariance of their arguments
obvious.
%
\section{Loop spaces of symmetric powers}
\label{loop-spaces-of-symmetric-Sec}
In this section, we will define symmetric powers of orbifolds and
study their loop spaces. In the case of a global quotient orbifold
$M\mmod G$, the $n^{th}$ symmetric power is again a global quotient
orbifold, namely $M^n\mmod (G\wr \Sn)$. We start by recalling some
well known facts about wreath products.
\subsection{Actions of wreath products}\label{wreath-Sec}
Let $G$ be a finite group.
Recall that the wreath product $G\wr\Sn$ has elements 
$$(\ug,\sigma)\in G^n\times\Sn$$
which compose as follows:
$$
  (g_1,\dots,g_n,\sigma)\cdot(h_1,\dots,h_n,\tau)=
  (g_1h_{\sigma^{-1}(1)},\dots,g_nh_{\sigma^{-1}(n)},\sigma\tau).
$$ 
Let $l$ be the order of $g_n\cdots g_1$ in $G$, and let $\sigma$ be the
$n$-cycle $(1\dots n)$. Then the order of $(\ug,\sigma)$ in $G\wr\Sn$
equals $ln$.

An element $(\uh,\tau)$ is in the centralizer of $(\ug,\sigma)$ if and only if 
$$
  \tau\sigma = \sigma\tau \quad\text{ and }\quad
  \forall i: g_{\sigma(\tau(i))}h_{\tau(i)} =  
  h_{\tau(\sigma(i))}g_{\sigma(i)}.
$$
Let $(\uh,\tau)$ be in the centralizer of $(\ug,\sigma)$, let
$(i_1,\dots,i_k)$ and $(j_1,\dots,j_k)$ be $k$-cycles in $\sigma$, 
and assume that $\tau(i_r) = j_{r+m}$.
Then
$$
  \forall r \in\ZZ/k\ZZ: g_{j_r}h_{j_{r-1}} = h_{j_r}g_{i_{r-m}}.
$$
In particular,
$$
  h_{j_k}g_{i_{k-m}}\cdots g_{i_{1-m}} = 
  g_{j_{k}}\cdots g_{j_{1}} h_{j_k},
$$
and the other $h_{j_r}$ are determined by $h_{j_k}$ and $\ug$. Note
that the element 
$$h_{j_k}g^{-1}_{i_{1-m}}\cdots g^{-1}_{i_{k}} = g_{j_1}^{-1}\cdots
g_{j_m}^{-1}h_{j_m}$$  
conjugates
$g_{j_k}\cdots g_{j_1}$ into $g_{i_k}\cdots g_{i_1}$.

Let $M$ be a right $G$-manifold. Then $G\wr\Sn$ acts on $\ux\in M^n$ via
$$\ux\cdot (\ug,\sigma) = 
(x_{\sigma(1)}g_{\sigma(1)},\dots,x_{\sigma(n)}g_{\sigma(n)}).$$
More generally, let ${\bf G}:=(G_0, G_1, s,t,u,i)$ be an orbifold
groupoid.
\begin{Def}
  We define the $n^{th}$ symmetric power of ${\bf G}$, denoted ${\bf
  G}\wr\Sn$, to be the groupoid ${\bf G}^n\rtimes\Sn$. Explicitely, 
  ${\bf G}\wr\Sn$ has objects $G_0^n$, and morphisms $G_1^n\times\Sn$. Its
  source sends the morphism $(\underline g,\sigma)$ to
  $(s(g_{\sigma(1)}),\dots,s(g_{\sigma(n)}))$, while its target sends it
  to $(t(g_1),\dots,t(g_n))$. The unit map sends the object
  $\underline x$ to $(u(x_1),\dots,u(x_n),1)$. The multiplication of
  $(\underline g,\sigma)$ and $(\underline h,\tau)$, where
  $g_{\sigma(i)} = h_i$, is given by
  $(g_1h_{\sigma^{-1}(1)},\dots,g_nh_{\sigma^{-1}(n)})$; and the
  inverse of $(\underline g,\sigma)$ is given by $(i(g_{\sigma(1)}),\dots,
  i(g_{\sigma(n)}),\sigma^{-1})$. 
\end{Def}
One checks that this definition is Morita-invariant.
We will be interested in the orbifold loop spaces of symmetric
powers. We start by considering the ``untwisted sector''.
\begin{Exa}[Symmetric powers of loop spaces]
  We have 
  $$
    \(\mL(M\mmod G)\)\wr\Sn = \(\coprod_{\underline{g}\in G^n} 
    \mL_{g_1}M\times\dots\times\mL_{g_n}M
    \)\mmod G\wr\Sn, 
  $$
  where $(\uh,\tau)\in G\wr\Sn$ acts on the right-hand side by
  $$
    (\gamma_1,\dots,\gamma_n)\mapsto
    (\gamma_{\sigma(1)}h_{\sigma(1)},\dots,\gamma_{\sigma(n)}h_{\sigma(n)}). 
  $$
\end{Exa}
  Assume now that $\mV$ is a loop bundle over $\mathcal LM$ with
  Fourier decomposition
  $$
    \mathcal V\at X \cong \sum_{j\in\ZZ}\(V_j\acted G\)q^\frac jN =:x 
  \quad
  $$  
%
Assume further that $\mV$ (and hence each $V_j$ is a right
$G$-bundle. Note that this is not the same setup as in Definition
\ref{S-bdl-Def}, here we are only considering the ``untwisted
sector'' 
$$
  \mL M\mmod G = \mL_1M\mmod C_1.
$$ 
The $n^{th}$ symmetric power $\mV^n\mmod \Sn$ is a $G\wr
\Sn$-equivariant bundle over $(\mL M)^n$, on which $\mathbb S^1$ acts
(diagonally). The Fourier decomposition of $(\mV\wr\Sn)\at{M^n}$ is 
%
%
%
\begin{eqnarray}\label{dN-Eqn}
      d_n(x)& := &   \sum_{\underline{j}\in\ZZ^n} 
    \(V_{j_1}\times0\times\dots\times 0\)
    q^{\frac{j_1}{N}}\oplus \dots \oplus (0\times\dots\times0\times
    V_{j_n}) q^{\frac{j_n}{N}}\\ 
\notag    &\cong&\sum_{j\in\ZZ}\sum_{k=0}^nW_{j,k}q^\frac jN,
\end{eqnarray}
where $W_{j,k}$ consists of the $\(\begin{matrix}{}n\\k\end{matrix}\)$
summands in which $V_j$ turns up exactly $k$ times and the other
summands are zero (in particular, $W_{j,0}=0$).
The summands of \eqref{dN-Eqn} are $G^n$-representations in the
obvious way, and the symmetric group acts on $W_{j,k}$ by permutation
of these summands. 
%

\subsection{Loop spaces of symmetric
  powers}\label{loops-of-symmetric-powers-Sec} 
We are now going to describe the loop space of $\MGSn$,
$$
  \mathcal L(\MGSn) = 
  \coprod_{[\ug,\sigma]}\mathcal L_{(\ug,\sigma)}(M^n)\mmod C(\ug,\sigma).
$$
Let $[\sigma]$ correspond to the partition $n=\sum_kkN_k$, and assume
that for each cycle of $\sigma$, we have fixed a ``first'' element
$i_1$, and thus a representation as 
$$(i_1,\dots,i_k).$$
\begin{Thm}\label{Ls-Thm}
  The
  component of $\mathcal L(\MGSn)$ corresponding to $[\ug,\sigma]$ is  
  homeomorphic to
  \begin{equation}\label{Lgs-Eqn}
    \Lgs\cong\prod_k\prod_{(i_1,\dots,i_k)}
    {}_k\mL_{g_{i_k}\cdot\dots\cdot g_{i_1}}M,
  \end{equation}
  where the second product runs over all $k$-cycles in $\sigma$.
  Let $(\uh,\tau)$ be an element of $C(\ug,\sigma)$ with $\tau(i_1) =
  j_{1+m}$ as in Section \ref{wreath-Sec}, and define its action on
  $\gamma\in{}_k\mL_{g_{j_k}\cdots g_{j_1}}M$ by
  $$\(\gamma\cdot(\uh,\tau)\)(t) :=
  \gamma(m+t)h_{j_k}g^{-1}_{i_{1-m}}\cdots 
  g^{-1}_{i_{k}}.$$
  Then $\gamma\cdot(\uh,\tau)$ is a path in ${}_k\mL_{g_{i_k}\cdots
    g_{i_1}}M$, and the homeomorphism is $C(\ug,\sigma)$-equivariant with
  respect to this action on the right-hand side. Moreover, it
  preserves the $\SS$-action by reparame\-tri\-za\-tion of paths.
  There is a canonical map from the target of (\ref{Lgs-Eqn}) with
  this action of $C(g,\sigma)$ to
  $$
    \prod_k\({}_k\mL(M\mmod G)\)\wr\Sigma_{N_k}. 
  $$
  Denote its composite with (\ref{Lgs-Eqn}) by $F_{(\ug,\sigma)}$. 
  As a map of orbifolds, $F_{(\ug,\sigma)}$
  is independent of the choices made.
\end{Thm}
  \begin{figure}[h]
    \centering
\includegraphics[scale=.25]{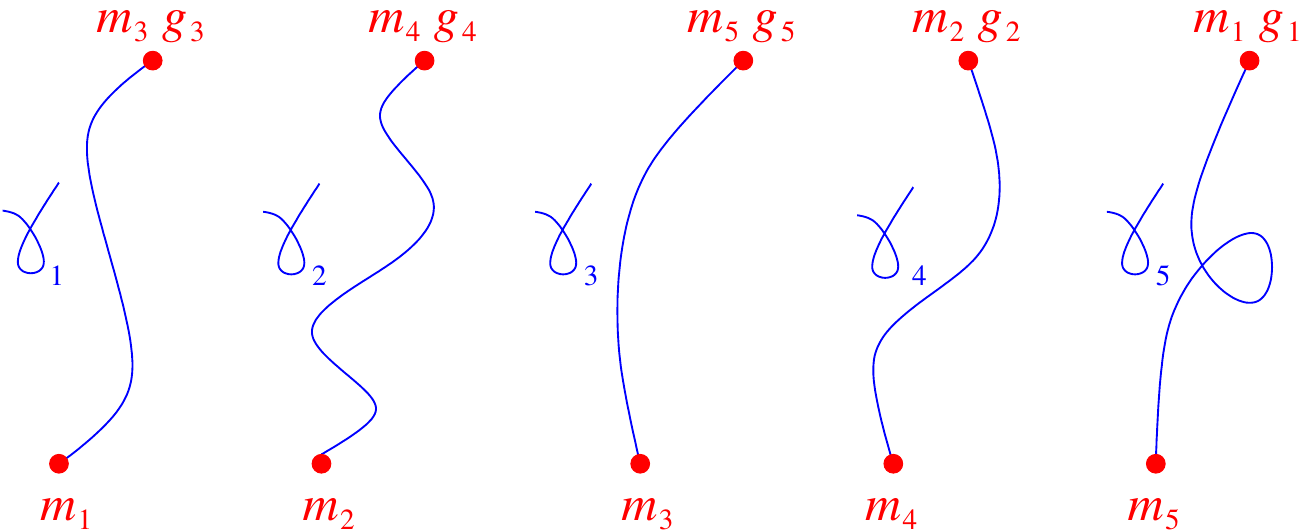}
\quad\quad\quad\quad\quad
\includegraphics[scale=.25]{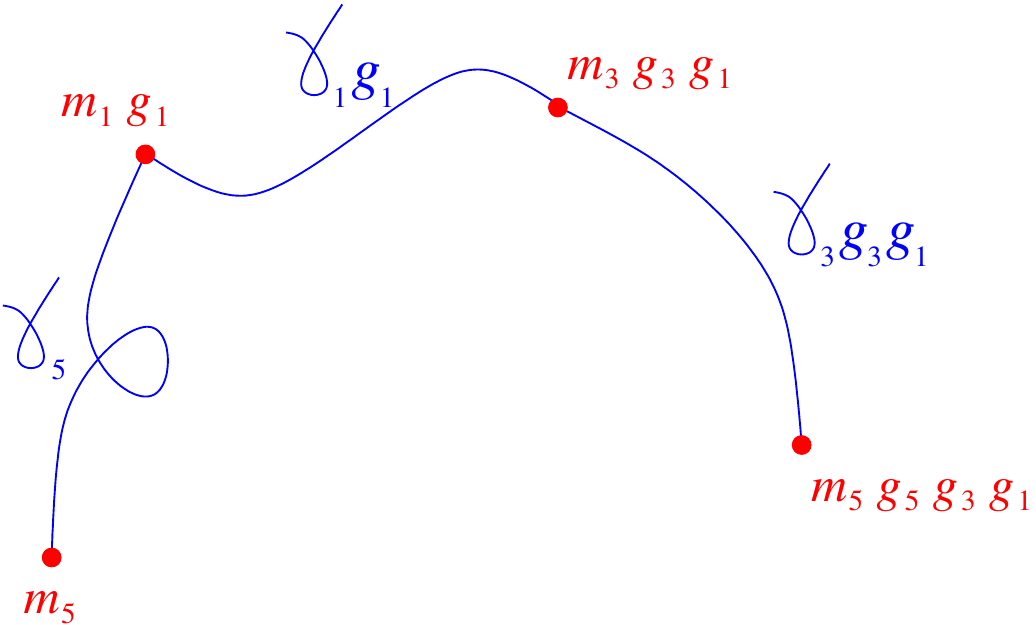}
\includegraphics[scale=.25]{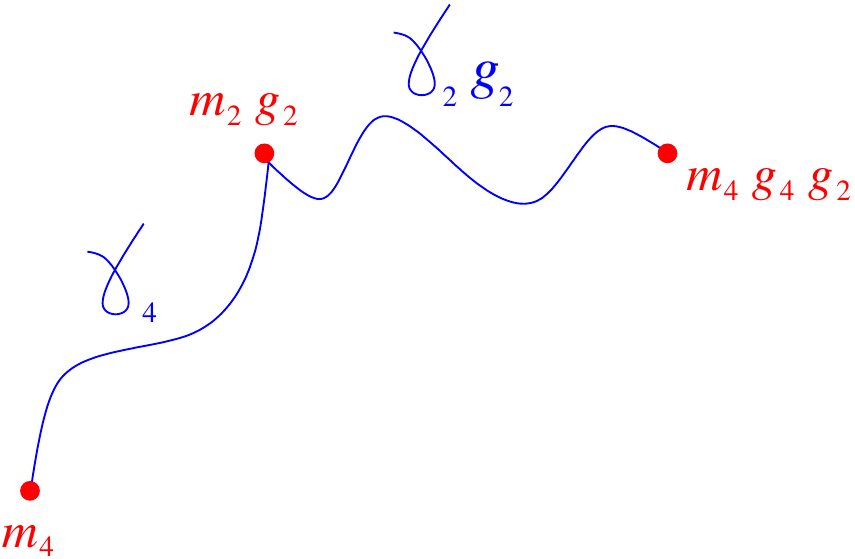}
    \caption{An element of $\mL_{((g_1,\dots, g_5);(135)(24))}M^5$ and its image in
      $\mL_{g_5g_3g_1}M\times\mL_{g_4g_2}M$.} 
\label{Loops-of-symmetric-powers:Fig}
  \end{figure}

\begin{Pf}{} 
The path
$$\gamma = (\gamma_1,\dots,\gamma_n)$$
in $M^n$ 
is in $\Lgs$ if and only if
$$\forall i: \gamma_i(1) = \gamma_{\sigma(i)}(0)g_{\sigma(i)}.$$

I.e., if $\gamma(0)=(m_1,\dots,m_n)$, then for $1\leq i\leq n$,
$\gamma_i$ is a path from $m_i$ to $m_{\sigma(i)}g_{\sigma(i)}$.
Let $(i_1,\dots,i_k)$ be a cycle in $\sigma$. The paths
$\gamma_{i_k}$, $\gamma_{i_1}g_{i_1}$ to 
$\gamma_{i_{k-1}}g_{i_{k-1}}\cdots g_{i_1}$ compose to a path
$$
  \gamma := \gamma_{i_k}*\gamma_{i_1}g_{i_1} *\dots
  *\gamma_{i_{k-1}}g_{i_{k-1}}\cdots g_{i_1}
$$
in ${}_k\mL_{g_{i_k}\cdots g_{i_1}}M$. This defines the homeomorphism
of the proposition (see Figure
\ref{Loops-of-symmetric-powers:Fig}).
Its inverse sends a path
$\gamma\in{}_k\mL_{g_{i_k}\cdots g_{i_1}}M$ to $(\gamma_{i_r})_{r=1}^k$ with 
$$\gamma_{i_k} = \gamma\at{[0,1]}$$
and
$$
  \gamma_{i_r}(t) := \gamma(r+t)g_{i_1}^{-1}\cdots g_{i_r}^{-1},
$$
for $1\leq r<k$. Let $(\uh,\tau)$ be as in the proposition.
Then the $i_r^{th}$ path of $(\gamma_i)_{i=1}^n\cdot(\uh,\tau)$ is
$\gamma_{j_{r+m}}h_{j_{r+m}}$,
and
$$
  \gamma_{j_{r+m}}h_{j_{r+m}}g_{i_r}\cdots g_{i_1} = 
    \gamma_{j_{r+m}}g_{j_{r+m}}\cdots g_{j_1} \( g_{j_1}^{-1}\cdots
    g_{j_m}^{-1}h_{j_m}\). 
$$
It follows that the centralizer acts as claimed.
Recall that 
$${}_k\mL(M\mmod G)\cong\coprod_{[g]}{}_k\mL_gM\mmod C^k_g$$
and that if $h$ is conjugate to $g$ the groupoids ${}_k\mL_gM\mmod C^k_g$
and ${}_k\mL_hM\mmod C^k_h$ are related by an isomorphism which is
canonical up to a natural transformation (i.e., it becomes canonical
in the category of orbifolds). To be specific, fix representatives of
the conjugacy classes of $G$, fix $k$, and fix an ordering
$\ui_1,\dots,\ui_{N_k}$ of the $k$-cycles of $\sigma$.
For $1\leq\alpha\leq N_k$, set $g_\alpha :=
g_{i_{\alpha,k}}\cdots g_{i_{\alpha,1}}$, and let $r_\alpha$ be the
representative of $g_\alpha$. Then on the $\alpha^{th}$ factor, the
canonical map of the proposition is represented by multiplication with
a group element $s_\alpha$ which conjugates $g_\alpha$ into $r_\alpha$
(and hence maps ${}_k\mL_{g_\alpha}M$ to ${}_k\mL_{r_\alpha}M$).
Let now $(\uh,\tau)\in C_G(\ug,\sigma)$ be as in the
proposition. Then $\tau$ defines a permutation $\rho\in\Sigma_{N_k}$
of the set of $k$-cycles of $\sigma$. Let $r$ be the representative of
$[g_{j_k} \cdots g_{j_1}]$. Since 
$h_{j_k}\cdot g_{i_{1-m}}\inv\cdots g_{i_{k}}\inv$ 
conjugates $g_{j_k} \cdots g_{j_1}$ into $g_{i_k}
\cdots g_{i_1}$, $r$ is also the representative of the latter, and
under our canonical map, right multiplication with $h_{j_k}\cdot
g_{i_{1-m}}\inv\cdots g_{i_{k}}\inv$  translates into multiplication
with an element of $C_r$ combined with the permutation $\rho$ of
factors of the target
$$\prod_{\alpha}{}_k\mL_{r_\alpha}M.$$
Rescaling a path by $t\mapsto t+m$ is the action of $m\in C_{r_\alpha}^k$.

Assume now that we have chosen a different first element, say  $i_r$,
of the cycle $(i_1,\dots,i_k)$. This leads to a different groupoid map
$F'_{(\ug,\sigma)}$. The factor of the target of $F_{(\ug,\sigma)}$
corresponding to 
this cycle is $\mL_{g_{i_k}\cdots g_{i_1}}M$, while that of
$F'_{(\ug,\sigma)}$ is $\mL_{g_{i_{r-1}}\cdots g_{i_r}}M$. There is a
natural isomorphism $I$ between the two maps, sending
$\gamma\in\mL_{g_{i_k}\cdots g_{i_1}}M$ to
$$
  I(\gamma)(t) = \gamma(t+r)g_{i_1}\inv\cdots g_{i_{r-1}}\inv.
$$
Hence $F$ and $F'$ define the same map of orbifolds.
Similarly, changing the order of the $k$-cycles of $\sigma$ translates
into a permutation in $\Sigma_{N_k}$.
%
%
%
\end{Pf}

The orbifold loop space $\mathcal L(\MGSn)$ can be viewed as the space of
$n$ strings moving in $M\mmod G$, their order does not matter, they are
either closed (in $M\mmod G$) or joining together to form longer closed strings.
\begin{Exa}\label{gsigma-Exa}
  Let $(\uh,\tau)$ be equal to $(\ug,\sigma)\in
  C_{(\ug,\sigma)}$. 
  Since $\sigma(i_1)=i_2$, and $g_{i_k}\cdot g_{i_k}\inv=1$,
  its action on $\mL_{g_{i_k}\cdots g_{i_1}}M$ 
  rotates the path $\gamma(t)$ to $\gamma(t+1)$. Hence $(\ug,\sigma)$ maps to
  $1\in C_r^k$ in the ${}_k\mL_rM$-factor of ${}_k\mL(M\mmod G)$ (see
  Figure \ref{Rotation-by-one}).

  \begin{figure}[h]
    \centering
$$
\begin{matrix}
  \raisebox{.7cm}{\includegraphics[scale=0.35]{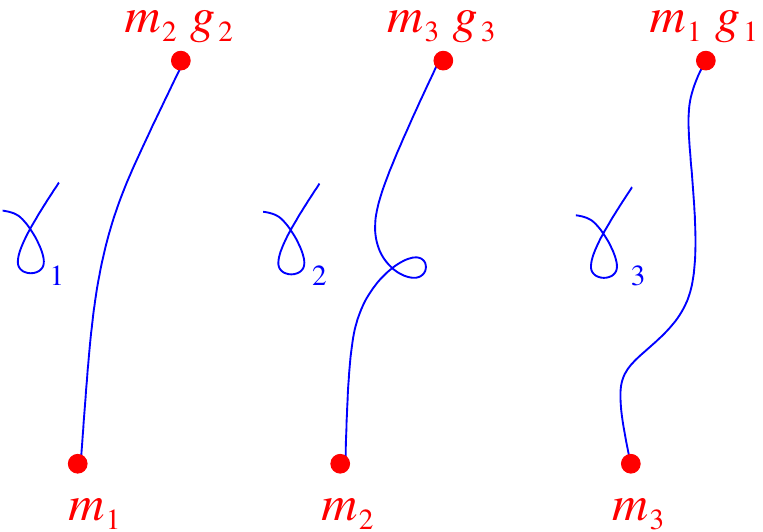}\quad}&
  \includegraphics[scale=0.35]{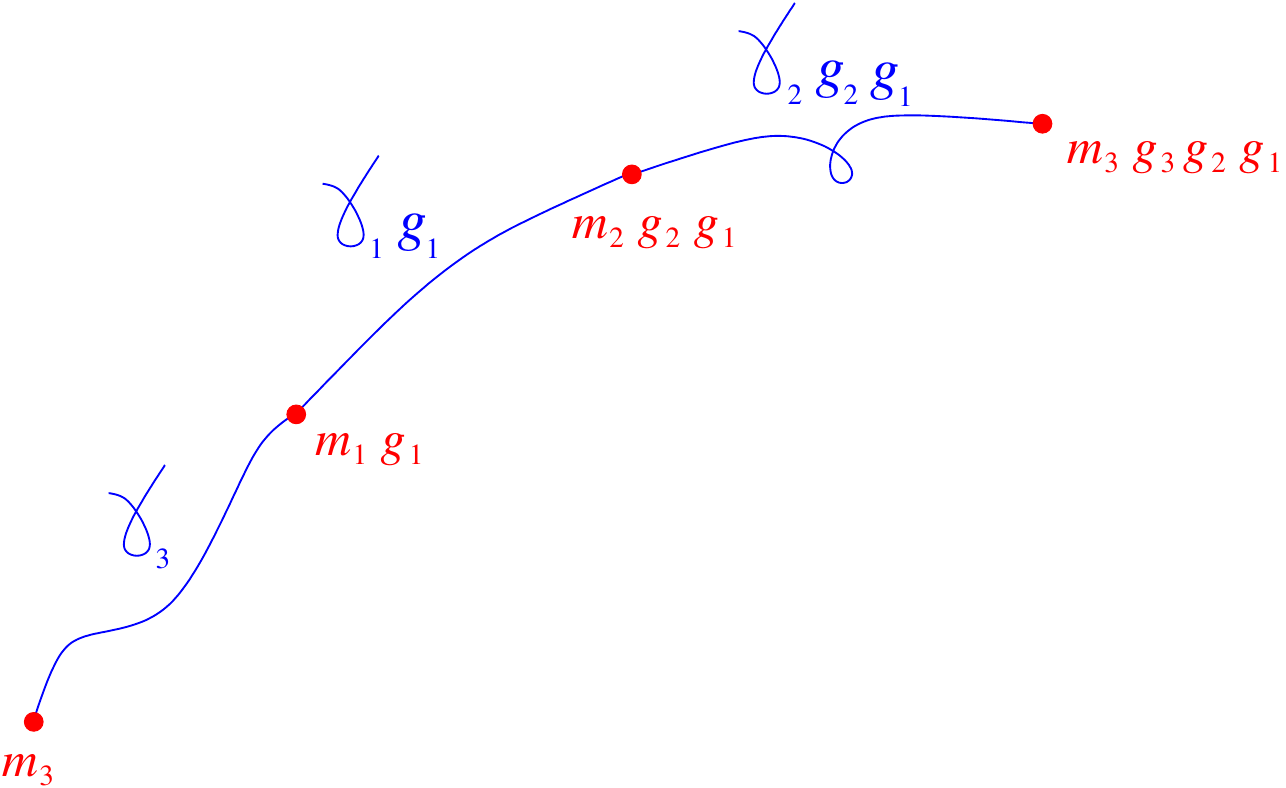}\\
  {\includegraphics[scale=0.35]{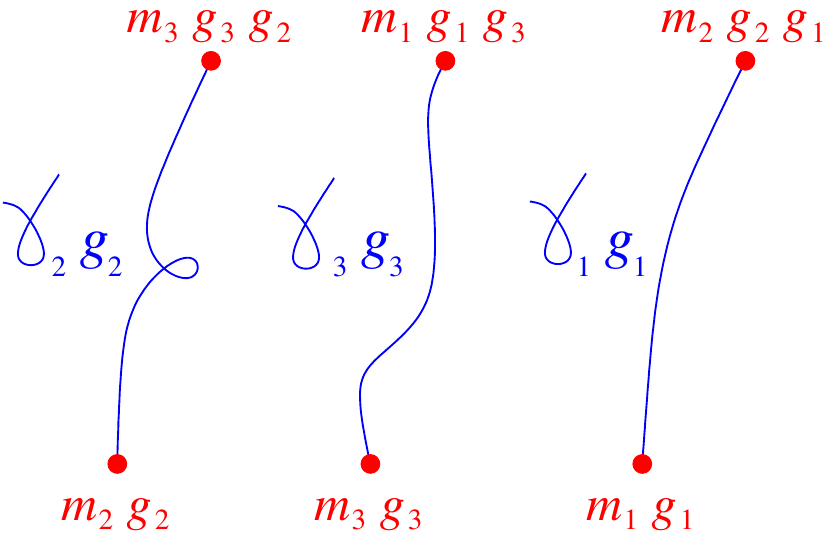}\quad}&
  \includegraphics[scale=0.35]{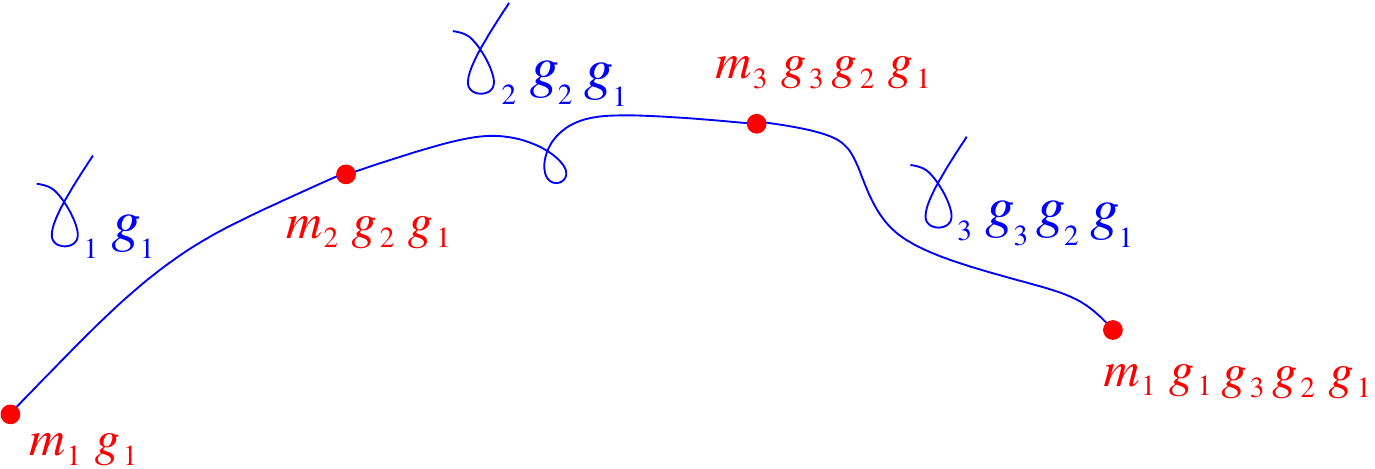}
\end{matrix}
$$
    \caption{The images in $\mL_{g_3g_2g_1}M$ 
      of an element of $\mL_{(g_1,g_2,g_3;(123))}M^3$ and its counterpart 
      under the action of $(g_1,g_2,g_3;(123))$
      differ by rotation by $1$.} 
\label{Rotation-by-one}
  \end{figure}
\end{Exa}

In the case that $G$ is the trivial group
the discussion in \cite{DMVV} yields a formula
 summarizing the above description of $\mathcal L(M\wr\Sn)$ for all $n$: 
write
$$
  \Sym^\bullet(tM) := \coprod_{n\geq 0}\MSn t^n
$$
for the ``total symmetric power'' of $M$. Here $t$ is a dummy variable.
\begin{Cor}
We have
$$
  \mathcal L\(\Sym^\bullet(tM)\) \cong
  \prod_{k\geq1}\Sym^\bullet((\kLM\rtimes(\ZZ/k))t^k),  
$$
and the inertia orbifolds are
$$
  \Lambda(\Sym^\bullet(tM)) := \coprod_{n\geq 0}\Lambda(M\wr\Sn)t^n 
  \cong\prod_{k\geq 1}\Sym^\bullet(t^kM\rtimes(\ZZ/k)),
$$
where the action of the groups $\ZZ/k\ZZ$ on $M$ is trivial.
\end{Cor} 
For non-trivial $G$, these product formulas are proved in
\cite{Tamanoi:Infinite}. 
\subsection{Fourier decompositions of loop bundles of symmetric products}
\label{Fourier-loop-symmetric-Sec}
In this section we will apply the construction of the previous section to
loop bundles and calculate its effect on Fourier decompositions.
The following two, slightly technical, definitions are motivated by
Example \ref{motivate-Dn-Exa} below. 
\begin{Def}\label{kx-Def}
  Let $\mV$ be an orbifold $\SS$-equivariant vector bundle over ${}_r\mL
  M\mmod G$. We write ${}_k\mV$ for the orbifold $\SS$-equivariant
  vector bundle over ${kr}\mL M\mmod G$ obtained by rescaling the
  $\RR/{kl}\ZZ$ action on $\mV_g$ with $\frac1k$ to obtain an action
  of $\RR/lkr\ZZ$. 
\end{Def}
Let $\mathcal V$ be as in the previous definition, and assume that its
decomposition over the constant $g$-loops is
$$
  \mathcal V\at {M^g}\cong V_0\bigoplus_{j\geq 1}V_j q^\frac j{rl}.
$$
Then the decomposition of ${}_k\mathcal V$ over $M^g$ is 
$$
  {}_k\mathcal V\at{M^g} \cong V_0\bigoplus_{j\geq1}V_j q^\frac j{krl},
$$
with the same $V_j$ now viewed as $C_g^{\RR/kr}$ representations,
where $C_g$ acts as before, and $\RR/klr\ZZ$ acts with rotation number
$j$. Hence the action of $g$ on $V_j$ equals that of
$kr\in\RR/klr\ZZ$, namely complex multiplication with $\zeta_l^j$.
\begin{Def}\label{D_n-Def}
  Let $\mV$ be an orbifold $\SS$-equivariant vector bundle over
  $\LMG$.  
  We define the $\SS$-eqivariant vector bundle $D_n(\mathcal V)$ over 
  $\mL((M\mmod G)\wr\Sn)$ by
  $$
    \(D_n(\mathcal V)\)_{(\ug,\sigma)} := F_{(\ug,\sigma)}^*\(
    \prod_k \({}_k\mV\)\wr\Sigma_{N_k}\),
  $$
  where $F_{(\ug,\sigma)}$ is the homoemorphism of Theorem \ref{Ls-Thm}.
\end{Def}
\begin{Exa}\label{motivate-Dn-Exa}
In the case that $\mathcal V=\mathcal L(V)$ is the loop space of a
bundle $V$ over 
$M\mmod G$, we have
$$D_n(\mV) \cong \mL((V\mmod G)\wr\Sn).$$
\end{Exa}

\label{fgs-page}
Let $(\ug,\sigma)$ be as above, and let
$f_{(\ug,\sigma)}=F_{(\ug,\sigma)}^{\SS}$ be the
restriction of $F_{(\ug,\sigma)}$ to the constant loops.
\begin{Prop}\label{Dn-Fourier-Prop}
  Let the orbifold $\SS$-equivariant vector bundle $\mV$ over
  $\mL(M\mmod G)$ have Fourier decomposition
  $$x:=\mathcal V\at{\Lambda(M)}\cong\bigoplus_{j\in\ZZ}(V_j\acted
  G)q^\frac j{|G|}$$ over the $G$-space $\coprod_gM^g$.
  Then the restriction of $D_n(\mV)$ to $(M^n)^{(\ug,\sigma)}$ has
  decomposition 
  $$f^*_{(\ug,\sigma)}\(d_{N_1}({}_1(x))\times
  d_{N_2}({}_2(x))\times\dots\times 
  d_{N_k}({}_k(x))\dots\)$$
  over the ${C_{G\wr\Sn}(\ug,\sigma)}$-space $\(M^n\)^{(\ug,\sigma)}$.
  Here $d_N$ is as in (\ref{dN-Eqn}), and ${}_k(x)$ is as in
  Definition \ref{kx-Def}.   
\end{Prop}
\begin{Pf}{}
  By the definition of $f_{(\ug,\sigma)}$, the diagram
  $$
    \xymatrix{{(M^n)^{(\ug,\sigma)}}\ar[r]^{f\phantom{xxxxxxx}}\ar[d]_\iota &
    {\prod\(\Lambda(M\mmod
    G)\)\wr\Sigma_{N_k}}\ar[d]^{\prod \iota_k} \\
    {\mL_{(\ug,\sigma)}M^n}\ar[r]^{F\phantom{xxxxxxx}} &
    {\prod\({}_k\mL(M\mmod G)\)\wr\Sigma_{N_k}}
  }
  $$
  commutes, where the products are over $k\in \mathbb N$, and the maps
  $\iota_k$ and
  $\iota$ are the fixed point inclusions. 
  Thus
  \begin{eqnarray*}
    (D_n\mV)\at{(M^n)^{(\ug,\sigma)}} & = &
    f_{(\ug,\sigma)}^*\(\prod_k
    \iota_k^*\(\({}_k\mV\)\wr\Sigma_{N_k}\)\)\\
    & = & f_{(\ug,\sigma)}^*\(\prod_k
    d_{N_k}\(\({}_k\mV\)\at{\Lambda(M\mmod G)}\)\)\\
    & = & f_{(\ug,\sigma)}^*\(\prod_k
    d_{N_k}\({}_kx\)\).
  \end{eqnarray*}
  Here $\prod$ stands for an external product, and $\Lambda(M\mmod G)$
  is the inertia orbifold as in Example \ref{translation-Exa}.
\end{Pf}
Proposition \ref{Dn-Fourier-Prop} motivates the following definition.
\begin{Def}\label{D_n-Def2}
  We define the map $D_n$ in $K_\Dev$ by
  \begin{eqnarray*}
    D_n\negmedspace : K_{\Dev,G}(M) & \to & K_{\Dev,G\wr\Sn}(M^n)\\
    x&\mapsto& f^*_{(\ug,\sigma)}\(d_{N_1}({}_1(x))\times
  d_{N_2}({}_2(x))\times\dots\times 
  d_{N_k}({}_k(x))\dots\).
  \end{eqnarray*}
\end{Def}
\subsection{Iterated symmetric powers}
We will need to understand
how our maps behave under iterated symmetric powers.
Let $$\ts=\tts\in\Sigma_m\wr\Sn.$$ For $1\leq i\leq n$
and $1\leq j\leq m$, we set
$$
  \ts(i,j) := (\sigma(i),\tau_{\sigma(i)}(j)).
$$
The bijection
\begin{eqnarray*}
  \{(i,j)\mid 1\leq i\leq n,1\leq j\leq m\} &\longrightarrow & \{1,\dots,nm\} \\ 
  (i,j)&\longmapsto& i+(j-1)n \\  
\end{eqnarray*}
induces an inclusion $$\map\iota{\Sigma_m\wr\Sn}\Sigma_{nm}.$$
If $(i_r,j_r)$ is a cycle of $\ts$, then the $i_r$ form a cycle of
$\sigma$. Let $\sigma$ be the cycle $(1\dots n)$, and set 
$$\tau = \tau_n\cdots\tau_1.$$ Then $(j_1\dots j_k)$ is a $k$-cycle of $\tau$ 
if and only if 
$$
  \( (1,\tau_1(j_k)) (2,\tau_2\tau_1(j_k)\dots(n,j_1)
   (1,\tau_1(j_1))\dots (n,j_2)(1,\tau_1(j_2))\dots\dots (n,j_k)
  \)
$$
is an $nk$-cycle of $\ts$.
\begin{Prop}\label{iterated-Prop}
Let 
$$(g_{i,j},\tau_i,\sigma)_{i,j}\in (G\wr\Sigma_m)\wr\Sn.$$
Then the following diagram commutes:
$$
\xymatrix{
{\mL_{(g_{i,j},\tau_i,\sigma)_{i,j}}(M^m)^n} \ar[r]\ar@{=}[d] &
{\prod\limits_l\prod\limits_{(i_1,\dots,i_l)}{}_l\mL_{(g_j,\tau)_j}M^m} \ar[d]\\
{\mL_{(g_{i,j},\iota(\tau_i,\sigma))_{i,j}}M^{mn}} \ar[r] &
{\prod\limits_{l,(i_1,\dots,i_l)}\prod\limits_{k,(j_1,\dots,j_k)}
{}_{kl}\mL_{g_{j_k}\cdots g_{j_1}}M} \\
}
$$
Here $(i_1,\dots,i_l)$ runs over all $l$-cycles of $\sigma$,
$$
  (g_j,\tau)_j := (g_{i_l,j},\tau_{i_l})\cdot\dots\cdot(g_{i_1,j},\tau_{i_1}),
$$
and the $(j_1,\dots,j_k)$ run over all $k$ cycles of $\tau$ (hence the
product in the lower right entry runs over all cycles of
$(\underline\tau,\sigma)$). The upper map is
$F_{(\underline{(\ug,\tau)},\sigma)}$, the right vertical map is
$\prod_{l,\underline i}F_{(\ug,\tau_i)}$, the left is $\mL\iota$ and
  the bottom is $F_{(\ug,\iota(\tau,\sigma))}$. (These are the maps
  defined in Theorem \ref{Ls-Thm}.)
\end{Prop}

\begin{figure}[h]
  \centering
  \includegraphics[scale=.27]{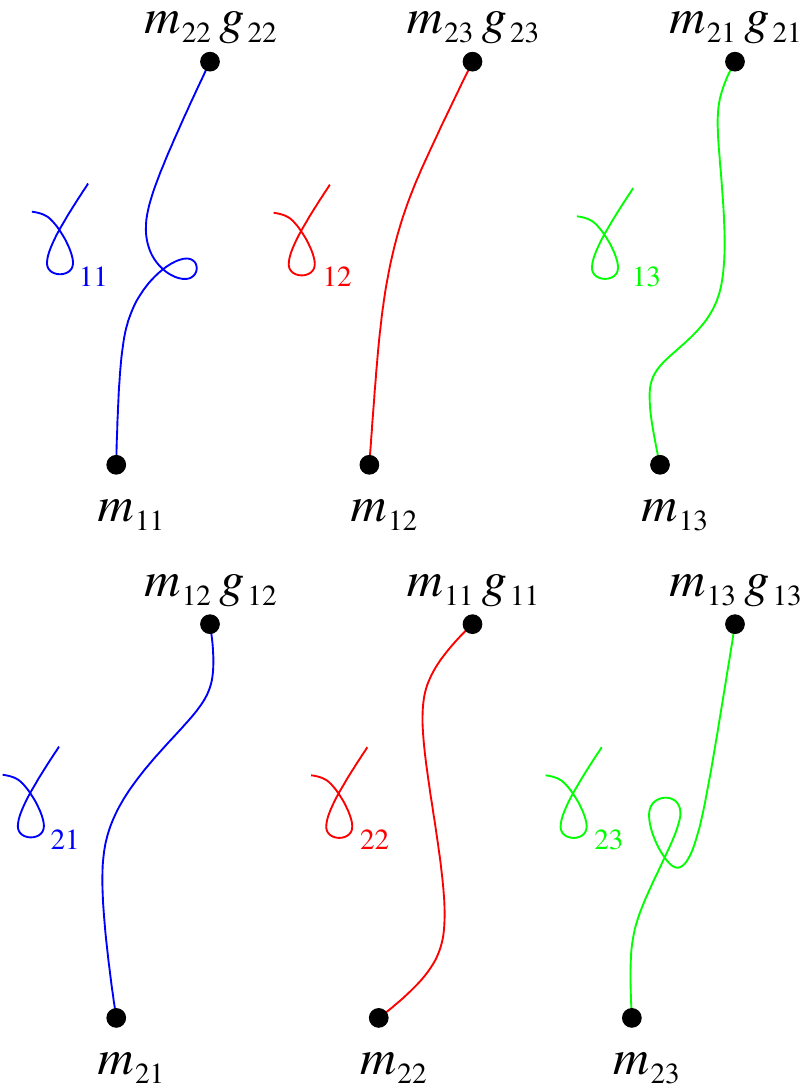}\quad\quad
  \includegraphics[scale=.27]{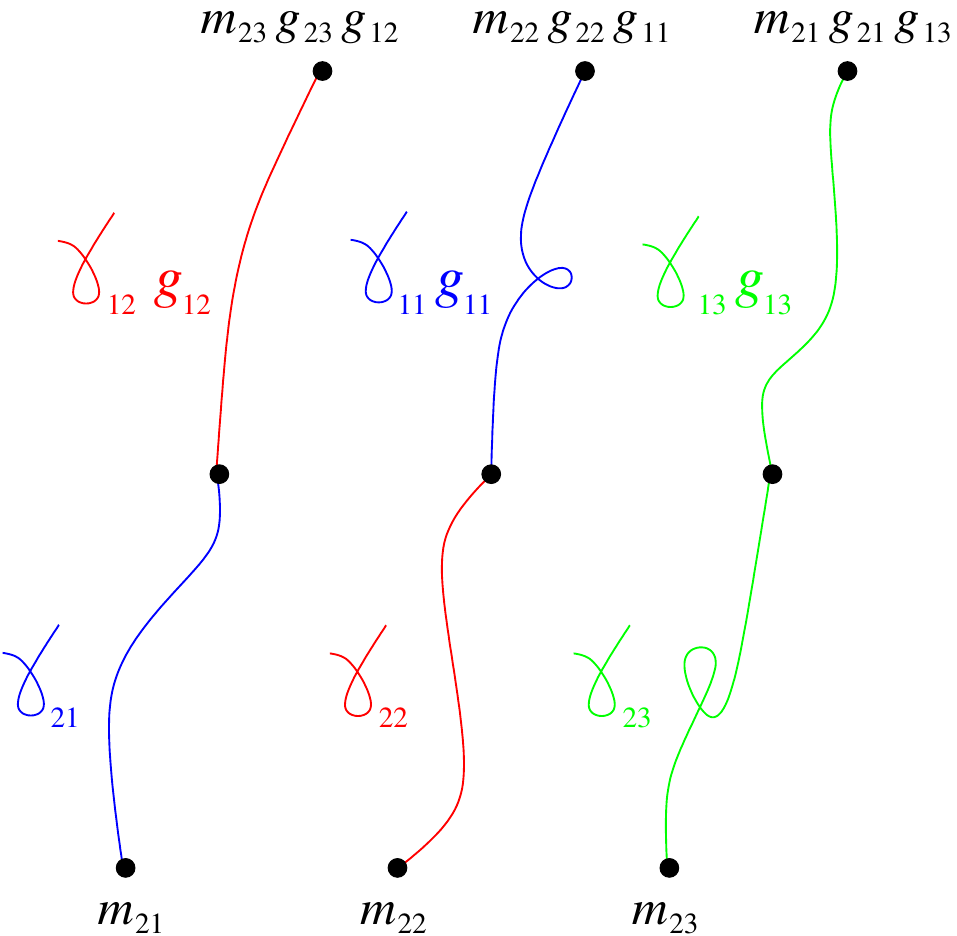}\quad\quad
 \includegraphics[scale=.27]{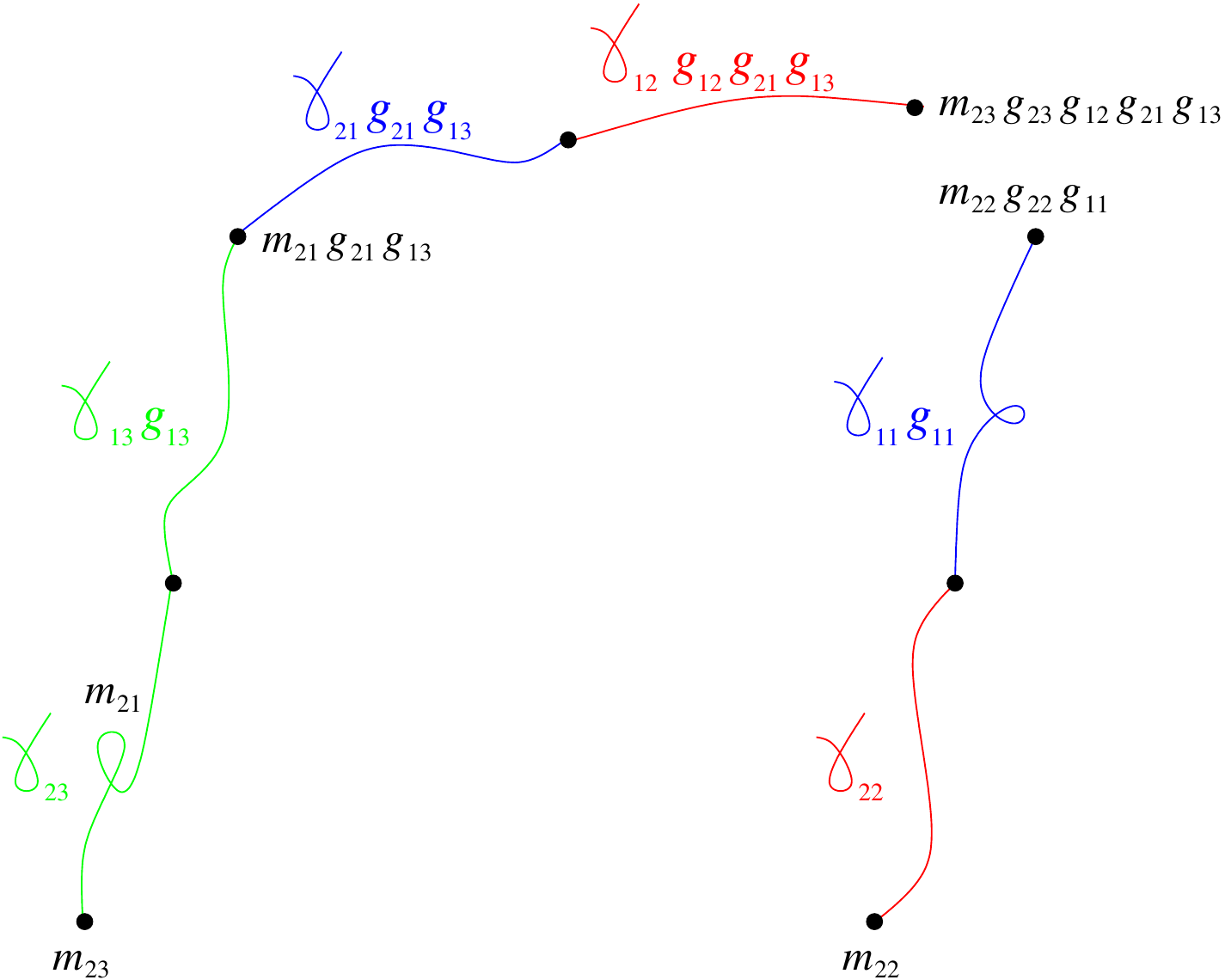}
  \caption{An example with $\sigma=(12)$, $\tau_1 =
    (12)(3)$, and $\tau_2=(123)$. Hence $\tau = (13)(2)$ and $g_1 =
    g_{21}g_{13}$ and $g_2 = 
    g_{22}g_{11}$ and $g_3=g_{23}g_{12}$.} 
  \label{twobythree}
\end{figure}

\begin{Pf}{} Without loss of generality,
  we may assume $\sigma = (1\dots n)$. 
  Then $F_{(\underline{(\ug,\tau)},\sigma)}$ sends
  $(\gamma_{i,j})_{i,j}$ to $(\gamma_j)_{j=1}^m$ with
  $$
    \gamma_j = \gamma_{n,j}*\gamma_{1,\tau_1(j)}g_{1,\tau_1(j)}
    *\dots*
    \gamma_{n-1,\tau_{n-1}\cdots\tau_1(j)}g_{n-1,\tau_{n-1}\cdots\tau_1(j)}\cdots 
    g_{1,\tau_1(j)}
  $$
  (see Figure \ref{twobythree}).
  Let 
  $$
    (g_1,\dots,g_m,\tau) := (g_{n,1},\dots,g_{n,m},\tau_n)\cdots
    (g_{1,1},\dots,g_{1,m},\tau_1).  
  $$
  Then $\tau = \tau_n\cdots\tau_1$, and
  $$g_j = g_{n,j}g_{n-1,\tau_n^{-1}(j)}\cdots
  g_{1,\tau_2^{-1}\cdots\tau_n^{-1}(j)}, $$
  and $\gamma_j$ is a path from $x_{n,j}$ to
  $x_{n,\tau(j)}g_{\tau(j)}$, where the $x_{i,j}$ are the starting
  points of the paths $\gamma_{i,j}$.
  Hence $(\gamma_j)_{j=1}^m$ is an element of  
  $$
    {}_n\mL_{(g_1,\dots,g_m,\tau)}M^m.
  $$
  Let now $(j_1,\dots,j_k)$ be a cycle of $\tau$. Then the component
  of $f_{m,G}$ corresponding to this cycle sends $(\gamma_j)$ to
  $$
    \gamma = 
    \gamma_{j_k}*\gamma_{j_{1}}g_{j_{1}}*\cdots*
    \gamma_{j_{k-1}}g_{j_{k-1}}\cdots g_{j_1}.
  $$
  For $r\in\ZZ/k\ZZ$, we have
  \begin{eqnarray*}
    \gamma\at{[rn+i,rn+i+1]} & = & \gamma_{j_r}g_{j_r}\cdots
    g_{j_1}\at{[i,i+1]}.
  \end{eqnarray*}
  Further, 
  \begin{eqnarray*}
    \gamma_{j_r}\at{[i,i+1]} & = &
    \gamma_{i,\tau_i\cdots\tau_1(j_r)}g_{i,\tau_i\cdots\tau_1(j_r)}
    g_{i-1,\tau_{i-1}\cdots\tau_1(j_r)}
    \cdots 
    g_{\tau_1(j_r)} \\  
    &  = & \gamma_{\ts^{rn+i}(n,j_k)}g_{\ts^{rn+i-1}(n,j_k)}
    \cdots g_{\ts^{rn+1}(n,j_k)},
  \end{eqnarray*}
  and
  $$
    g_{j_s} = g_{\ts^{sn}(n,j_k)}\cdots g_{\ts^{(s-1)n+1}(n,j_k)}.
  $$
\end{Pf}
\begin{Cor}
  Let $\iota$ be the inclusion of $\Sigma_m\wr\Sn$ in $\Sigma_{mn}$
  described above. Then we have
  $$\res{}\iota \circ D_{nm} \cong D_n\circ D_m.$$
\end{Cor}
%
%
\section{Power operations}\label{power-operations-Sec}
The goal of this section is to define power operations on Devoto's
equivariant Tate $K$-theory. We want to define them in such a way that
the equivariant Witten genus becomes an $H_\infty^2$ map. In other
words, we need our power operations to take Thom classes
$u_{\AA}^{\str}(V)$ to the Thom classes $u_{\AA}^{\str}(V\wr \Sn)$. It
turns out that the discussion of the previous section dictates our
definitions. In our treatment of power operations, we will once again
proceed in two steps. First, we 
define ``{\em Atiyah power operations}'' on the untwisted sector
$K_G(X)\ps q$. Then we use the Atiyah power operations to define
``{\em stringy power operations}'' on all of $K_{\Dev, G}(X)$.
We briefly recall the general setup of equivariant power
operations from \cite[Def.4.3]{Ganter:thesis}. 
Let 
$\{E_G\mid G \text{ finite}\}$ 
be a compatible family of equivariant cohomology theories in the sense of
\cite[II.8.5]{Lewis:May:Steinberger}, and write $\EG X$ for $E_G^0(X)$. 
We also ask that our family has unitary, commutative and associative
external products $\boxtimes$ that are natural in (stable) maps of $X$ and $Y$. 
%
%
\begin{Def}\label{hinf-Def}
  An $\hinf$-structure on $E$ is given by a collection of natural maps
  $$
    \map{P_n}{E_G(X)}{E_\GS(X^n)}
  $$ 
  called {\em power operations} satisfying the following conditions: 
  \begin{enumerate}
    \renewcommand{\labelenumi}{(\alph{enumi})}
    \item $P_1 = \id$ and $P_0(x) = 1$,
    \item the (external) product of two power operations is 
      $$P_j(x)\boxtimes P_k(x) = 
        \res{\Sigma_{j+k}}{\Sigma_j\times\Sigma_k} (P_{j+k}(x)),
      $$
    \item the composition of two power operations is
      $$
        P_j(P_k(x)) = \res{\Sigma_{jk}}{\Sigma_k\wr\Sigma_j}\(P_{jk}(x)\),  
      $$
    \item and the $P_j$'s preserve (external) products:
    $$P_j(x\boxtimes y) = \res{\Sj\times\Sj}\Sj(P_j(x)\boxtimes P_j(y)),$$
  \end{enumerate} 
  where the restriction is along the map
  $$
    \left[((X\acted G)^2)\wr\Sigma_j\right] \longrightarrow \left[(X\acted
    G)\wr(\Sj\times\Sj)\right]\cong
    \left[((X\acted G)\wr\Sj)^2\right]. 
  $$
\end{Def}
We also recall the graded ring
$$
  S_E(X) = \cOplus E_{G\wr\Sn}(X^n)t^n.
$$
\label{S_E-page}
Here $t$ is a dummy variable keeping track of the grading, and the
multiplication of two elements of degree $m$ and $n$ respectively, is
given by the external product, followed by the transfer
$$
  \ind{\Sigma_m\times\Sigma_n}{\Sigma_{m+n}}.
$$
The ring $S_E(X)$ is the target of the total power operation 
$$
  P = \sum_{n\geq 0}P_n t^n.
$$
\subsection{Atiyah power operations}
In this section we work over the untwisted sector, i.e., with the
compatible family
$$
  E_G(X) = K_G(X)\ps q.
$$
On $K_G(X)$ we have Atiyah's power operations
\begin{eqnarray*}
  P_n\negmedspace : K_G(X) &   \to   & K_{G\wr \Sn}(X^n)\\
                       {[V]} &\mapsto  & {[V^{\boxtimes n}]}.
\end{eqnarray*}
Here $\Sn$ acts with a sign on $V^{\boxtimes n}$.
\begin{Def}\label{Atiyah-powerop-Def}
  We define the total Atiyah power operation on $K_\Tate$ by
  \begin{eqnarray*}
    \Ptop\negmedspace :K_G(X)\ps q & \to & \cOplus K_{G\wr\Sn}(X^n)\ps qt^n\\
                   x_jq^j       & \mapsto & \sum_{n\geq 0}
                                     x_j^{\tensor n}q^{jn} t^n, \text{and} \\
            \sum_{j\geq 0}x_jq^j & \mapsto &  \prod_{j\geq 0}\Ptop(x_jq^j).
  \end{eqnarray*}
\end{Def}
The coefficient of $t^n$ is called the $n^{th}$ Atiyah power operation, it 
is given by
$$
  P_n^{\on{top}}(x) = \sum_{\ui}P_\ui(x),
$$
where $\ui = (i_j)_{j\geq 0}$ runs over all sequences of natural
numbers whose elements 
add up to $n$, and
$$
  P_\ui(x) = \prod_jP_{i_j}(x_j)q^{i_jj}.
$$
Note that for each $\ui$, this product has only finitely many
non-trivial factors. It is taken in the graded ring $S_{K_G}(X)$ and
hence involves 
the induced representations
$$
  \ind\ui n :=\ind{\Sigma_{i_1}\times\dots\times\Sigma_{i_j}\times \cdots}{\Sn}
$$
The transfer $\ind\ui n$ adds all possible shuffles of the partition
$\ui$. More precisely, we have
\begin{equation}
  \label{P_n-Eqn}
  P_n^{\on{top}}(x) = \sum_{\uj\in\mathbb N^n} x_{j_1}q^{j_1}
  \boxtimes\dots\boxtimes x_{j_n}q^{j_n} ,
\end{equation}
and $\Sn$ acts (with a sign) by permutation of the summands and the
factors therein.
The infinite sum $P_n^{\on{top}}$ is a well defined formal power
series in $q$. 
\begin{Prop}
  The $P^{\on{top}}_n$ satisfy the axioms of Definition \ref{hinf-Def}.
\end{Prop}
\begin{Pf}{}
  Part (a) is clear. Part (b) and (c) follow from Equation
  (\ref{P_n-Eqn}).
  Let $x\in K_G(X)\ps q$ and $y\in K_H(Y)\ps q$. Then the external
  product of $x$ and $y$ is
  $$
    x\boxtimes y = (\sum_{i}x_iq^i)\boxtimes(\sum_j y_jq^j) =
    \sum_{i,j}(x_i\boxtimes y_j)q^{i+j},
  $$
  and
  $$
    P_n(x\boxtimes y) = \sum_{\ui,\uj}(x_{i_1}\boxtimes
    y_{j_1})q^{i_1+j_1}
    \boxtimes  
    \dots\boxtimes (x_{i_n}\boxtimes y_{j_n})q^{i_n+j_n},  
  $$
  Under the isomorphism $$(X\times Y)^n \cong X^n\times Y^n,$$
  this is identified with 
  $$
    \res{(G\wr\Sn)\times(H\wr\Sn)}{(G\times H)\wr\Sn}P_n(x)\boxtimes P_n(y).
  $$
\end{Pf}

\subsection{The topological Witten genus of orbifolds and its product
  formula}
\begin{Def}
  
\end{Def}
Let $V$ be a $G$-equivariant $\Spin$ (respectively complex) vector
bundle, and set 
$$
  u_G^\top(V) := 
  u^\AA_G(V)\cdot\prod_{j=1}^\infty
  \Lambda_{-q^j}V_\CC. 
$$
We write 
$$
  \map{\phi_G^\top}{\on{MSpinP}_G}{K_G\ps q}.
$$
(respectively with $\on{MSPinP}_G$ replaced by $\on{MUP}_G$)
for the corresponding genus (c.f.\ Proposition \ref{Okonek-Prop}).
\begin{Prop}\label{top-H_oo-Prop}
  The map $\phi^\top_G$ is $H_\infty^4$ ($H_\infty^2$ in the complex case).
\end{Prop}
\begin{Pf}{}
  By \cite[4.6]{Ganter:thesis},\cite[(A4)]{tomDieck}, it suffices to
  show that
  $$
    u^\top_{G\wr\Sn}(V^n) = P_n(u_G^\top(V))
  $$
  for any $\Spin(4k)$-bundle $V$.
  This follows from Axiom (d),
  $\dim V^n = 4nk$, the known analogue for $u^\AA$ and $u^\td$: 
  $$u_{G\wr\Sn}^\AA(V^n) =
  P_n(u_G^\AA(V))$$ \cite{BMMS},
  and the fact that
  $$\Lambda_{t}(V^n) \cong (\Lambda_tV)^{\boxtimes n},$$
  where $\Sigma_n$ acts by permuting the factors inside on the left-hand
  side, and by permuting the factors outside (with a sign) on the
  right-hand side.
\end{Pf}
\begin{Def}
  We define the topological Witten genus of a global quotient orbifold
  $M\mmod G$ by 
  $$
    \phi_\top(M\mmod G) := \GG G\sum_{g\in G}\tr{\phi_G^\top(M)},
  $$
  where $\phi$ stands for $\AA$ or $\td$.
\end{Def}
\begin{Def}
  We define the (topological) total symmetric power operation $\Stop$ on
  $K_\Tate$ to be the composite
  $$
    \Stop\negmedspace : K_G(X)\ps q\to\cOplus K_{G\wr\Sn}(X^n)\ps qt^n
    \to\cOplus K_{G\times\Sn}(X)\ps qt^n
    \to K_G(X)\ps q\ps t,
  $$
  where the first map is $\Ptop$, and the second map is restriction along
  the diagonals of $X^n$ and $G^n$. 
  The last map $\eps$ is (in degree $n$) given by
  $$
    \eps_n = \map{\id\tensor\langle-,1\rangle_\Sn}{K_G(X)\tensor R(\Sn)}{K_G(X)}.
  $$
\end{Def}
Note that the total power operation takes sums to products, and that
$\eps$ is a ring map. Hence $S_t$ is exponential.
In the special case that $X$ is a point and $G$ is the trivial group,
we obtain the following corollary of Proposition \ref{top-H_oo-Prop}:
\begin{Cor}[product formula for the topological Witten genus]
  If $M$ is a $\on{Spin}(4k)$ or complex manifold with Witten genus
  $\phi(M) = \sum_jc(j)q^j$, then 
  $$
    \sum_{n\geq 0}\phi_\top(M^n\mmod\Sn)t^n = \prod_{j}\(\frac1{1-tq^j}\)^{c(j)}.
  $$
\end{Cor}

\begin{Pf}{}
  Since $\phi^\top$ preserves power operations, we have
  $$\sum_{n\geq 0}\phi^\top_\Sn(M^n)t^n = P^\top(\phi(M)).$$ Applying
  $\eps$ to both sides, we get
  $$
    \sum_{n\geq 0}\phi_\top(M^n\mmod\Sn)t^n = \Stop(\phi(M)).
  $$   
  Since $\Stop$ is exponential, it suffices to observe that $\Stop(q^j) =
  \frac1{1-tq^j}$ to complete the proof.
\end{Pf}
\subsection{Stringy power operations}
%
The following definitions should be compared to those of Section
\ref{Fourier-loop-symmetric-Sec}.
\begin{Def}\label{kx-ring-Def}
  Let ${}_k(-)$ be the ring map
  \begin{eqnarray*}
    {}_k(-)\negmedspace : {K_{\Dev, G, r}}&\to&{K_{\Dev, G, rk}}\\
    \sum V_jq^\frac j{r|g|} &\mapsto& \sum V_jq^\frac j{rk|g|}. 
  \end{eqnarray*}  
  Here the action of $C_g$ on $V_j$
  remains unchanged, while $1\in\ZZ/rk|g|\ZZ$ acts as
  $\zeta_{rk|g|}^j$ on the $j^{th}$ coefficient of the $[g]^{th}$
  summand on the right-hand side.
%
\end{Def}
%
%
\label{atiyah-powerop-in-Kdev-page}
In the following, for $x\in K_{\Dev,G,k}(X)$, by $P^\top_n(x)$ we
will mean the image of 
$x$ under the composite 
$$
  K_{\Dev,G,k}(X)\to K_{G\times\ZZ/kL}(\coprod_gX^g)\ps{q^\frac1{kL}}\to
  K_{(G\times\ZZ/kL)\wr\Sn}((\coprod_gX^g)^n)\ps{q^\frac1{kL}}. 
$$
Here $L$ is the order of $G$, and the second map is $P_n^\top$.
\begin{Def}\label{stringy-P-Def}
  Let $x\in K_{\Dev,G}(X)$.
  We define $$P_n^\str(x)\in K_{Dev,G\wr\Sn}(X^n)$$ by
  $$
    \(P_n^\str(x)\)_{(\ug,\sigma)} := f^*_{(\ug,\sigma)}\(
    P^\top_{N_1}\({}_1(x)\)\boxtimes\dots\boxtimes
    P^\top_{N_k}\({}_k(x)\)\boxtimes 
    \dots
    \).
  $$
\end{Def}
We need to verify that $P_n^\str$ indeed takes values in $K_{Dev,G\wr\Sn}(X^n)$.
First, we have
\begin{eqnarray*}
  P^\top_{N_k}({}_k(x))& \in &
  K_{(G\times\ZZ/kL)\wr\Sigma_{N_k}}((\coprod_{h\in G}X^h)^{N_k}) \ps{q^\frac1{kL}}  \\ 
  &\cong &\bigoplus_{[\uh]}K_{\on{Stab}(\uh)}(\prod_{i=1}^{N_k}
  X^{h_i})\ps{q^\frac1{kL}}, 
\end{eqnarray*}
where the direct sum is over the orbits of the action of $(G\times
\ZZ/kL)\wr\Sigma_{N_k}$ on $G^{N_k}$.
In fact, the $\uh^{th}$ summand of $  P^\top_{N_k}({}_k(x))$ is a
power series in $q^\frac1{k|\uh|}$, since it is obtained by
multiplying power series in $q^\frac1{k|h_i|}$ for $1\leq i\leq N_k$.
Further, recall from Section \ref{loops-of-symmetric-powers-Sec} that 
$$
  (X^n)^{(\ug,\sigma)}\cong\prod_{k\geq1}\prod_{(i_1,\dots,i_k)}
  X^{g_{i_k}\cdots g_{i_1}},  
$$
where the second product runs over the $k$-cycles of $\sigma$. Hence,
for fixed $k$, the
target of the $k^{th}$ factor of $f_{(\ug,\sigma)}$ is
$\prod_{i=1}^{N_k} X^{h_i}$, where $i$ runs over the $k$-cycles of
$\sigma$, and $h_i=g_{i_k}\cdots g_{i_1}$.
By Example \ref{gsigma-Exa}, $(\ug,\sigma)$ acts by multiplication with
$\zeta_{k|h_i|}^j$ on the coefficient of $q^{\frac{j}{k|h_i|}}$ of
${}_kx\at{X^{h_i}}$.
Since 
$|(\ug,\sigma)|$ is the smallest common multiple of all 
$k\cdot |g_{i_k}\cdots g_{i_1}|$ such that $k\in\mathbb N$ and $(i_1\dots i_k)$
is a $k$-cycle of $\sigma$, it follows that
$\(P_n^\str(x)\)_{(\ug,\sigma)}$ is a power series in
$q^\frac1{|(\ug,\sigma)|}$, where $(\ug,\sigma)$ acts by
multiplication with $\zeta_{|(\ug,\sigma)|}^j$ on the coefficient of
 $q^{\frac{j}{|(\ug,\sigma)|}}$.
%
\begin{Thm}
  The $P^{\str}_n$ satisfy the axioms of Definition \ref{hinf-Def}.
\end{Thm}
%
\begin{Pf}{}
  Axiom (a) is clear. For Axiom (b), let
  $(\sigma,\tau)\in\Sigma_m\times\Sn$. Then the set of $k$-cycles of the
  element $(\sigma,\tau)\in\Sigma_{m+n}$ is identified with the
  disjoint union of the 
  sets of $k$-cycles of $\sigma$ and $\tau$, and thus
  \begin{eqnarray*}
    \(P^\str_{n+m}(x)\)_{(\ug,(\sigma,\tau))} &=& 
       \(P_m^\str\)_{(g_1,\dots,g_m,\sigma)} \boxtimes
       \(P_n^\str\)_{(g_{m+1},\dots,g_{m+n},\tau)}\\ 
     & = &  \(P_m^\str \boxtimes P_n^\str\)_{(\ug,(\sigma,\tau))}. 
  \end{eqnarray*}
  Axiom (c) follows from Proposition \ref{iterated-Prop}: With the
  notation of that Proposition, and omitting the equivariant
  information, we have
  \begin{eqnarray*}
    \(P^\str_n\(P^\str_m(x)\)\)_{(g_{i,j},\tau_i,\sigma)} &=&
    \prod_{l,(i_1,\dots,i_l)}{}_l\(\prod_{k,(j_1,\dots,j_k)}
    {}_kx_{g_{j_k}\cdots g_{j_1}}\) \\
    &=&   P^\str_n(x)_{(g_{i,j},\tau_i,\sigma)},
  \end{eqnarray*}
  Here all the products denote the external tensor product
  $\boxtimes$. For the non-equivariant information it makes no difference
  whether ${}_l(-)$ is inside or outside the product, and when
  restricting along the inclusion of centralizers
  $$
    C_{G\wr\Sm\wr\Sn}(g_{i,j},\tau_i,\sigma)\sub C_{G\wr\Sigma_{mn}}(g_{i,j},\tau_i,\sigma)
  $$
  we restrict that part along the inclusion
  $\ZZ/l\ZZ\sub(\ZZ/lk\ZZ)^{N_k}$.
  To prove Axiom (d), let $x\in K_{\Dev,G}(X)$ and $y\in
  K_{\Dev,H}(Y)$. Then $x\boxtimes y$ is the element of
  $K_{\Dev,G\times H}$ with 
  $$
    (x\boxtimes y)_{[g,h]} = x_{[g]}\boxtimes
    y_{[h]},
  $$ 
  and the statement follows from the analogous statement for
  the $P_{N_k}^\top$.   
\end{Pf}

%
%
In this section, we will work with the following variant of the
stringy equivariant Witten genus:
\begin{Def}\label{stringy-Thom-Def}
  We define
  $$
    \map{\td^\str_G}{\MU_G}{K_\Dev}
  $$
  to be the map of spectra corresponding to the Thom class 
  $$
    u_{\td,G}^\str(V)_{[g]}:=u_{C_g}^\td(V^g)\cdot 
    \prod_{j\geq 1}\Lambda_{-1}\((V_\CC)_{\zeta_l^j} q^\frac jl\).
  $$
\end{Def}
The computation of $\td^\str_{[g]}(M)(h)$ is the same as that for
$\AA^\str$ with the terms
$$
  q^{-\frac d{12}+\frac\ageg2}e^{i\pi\medspace\ageh}\prod e^{\tilde x_i}
$$
missing. For $\on{SU}$-bundles, we have $\sum \tilde x_i = 0$, so that
the missing factors are just a renormalization constant.
\begin{Thm}\label{H_oo-Thm}
  The map $\td^\str_G$ is 
  $H_\infty^2$.
\end{Thm}
\begin{Pf}{}
  We write $u_G^\str$ for $u_{\td,G}^\str$ and $u$ for $u^{\td}$.
  We need to prove
  $$
    u_{G\wr\Sn}^\str(V^n) = P_n^\str(u_G^\str(V)).
  $$
  Let $(\ug,\sigma)\in G\wr\Sn$. Then
  $$
    V^n\at{(X^n)^{(\ug,\sigma)}}\cong\prod_k\prod_{\ui}
    V^{\oplus k}\at{X^{g_{i_k}\dots g_{i_1}}},
  $$
  where the product signs denote external direct sums and the second
  one runs over the $k$-cycles of $\sigma$.
  On the $\ui^{th}$ summand of this, $(\ug,\sigma)$ acts by
  $$(v_1,\dots,v_k)\mapsto(v_2g_{i_2},
  \dots,v_kg_{i_k},v_1g_{i_1}).$$
  Hence $\zeta$ is an eigenvalue of this action with eigenvector
  $(v_1,\dots,v_k)$ if and only if $\zeta^k$ is an eigenvalue of
  $g_{i_k}\dots g_{i_1}$ with eigenvector $v_k$, and $v_{j} =
  v_{j+1}g_{i_{j+1}}$ for $1\leq j\leq k$. Thus 
  $$
    (V^{\oplus k}\at{X^{g_{i_k}\dots g_{i_1}}})_\zeta\cong(V\at{X^{g_{i_k}\dots g_{i_1}}})_{\zeta^k},
  $$
  and   
  $$
    \(V^n\at{(X^n)^{(\ug,\sigma)}}\)_\zeta\cong\prod_k\prod_{\ui}
    \(V\at{X^{g_{i_k}\dots g_{i_1}}}\)_{\zeta^k}.
  $$
  Now Definition \ref{stringy-Thom-Def} gives
  \begin{eqnarray*}
    u^\str_{(\ug,\sigma),\alpha}(V^n) &=& 
    u_{C_{(\ug,\sigma)}}(V_1)\cdot\prod_{j=1}^{\infty}
    \Lambda_{-q^{\frac{j}{|\ug,\sigma|}}}\(\(V_\CC\)_{\zeta_{|\ug,\sigma|}^j}\)\\ 
    &=& \prod_{k,\ui}\(
    u_{C_{g_i}^k}(V_1)\cdot\prod_{j=1}^{\infty}
    \Lambda_{-q^{\frac{j}{k|g_i|}}}\(\(V_\CC\)_{\zeta_{|g_i|}^j}\)\),\\ 
  \end{eqnarray*}
  where $g_i=g_{i_k}\cdots g_{i_1}$, in the first line, $V$ and
  $V_\CC$ get restricted to $(X^n)^{(\ug,\sigma)}$ and in the second
  line, they get restricted to $X^{g_{i_k}\cdots g_{i_1}}$.
\end{Pf}

\subsection{Hecke operators}\label{Hecke-Sec}
In this section we prove that the power operations $P^\str_n$ are
elliptic in the sense explained in the introduction.
We recall 
from \cite{Ganter:thesis} that any cohomology theory with
power operations and a level $2$ Hopkins-Kuhn-Ravenel theory has Hecke
operators, acting as (internal) cohomology operations:
\begin{Def}
The $n^{th}$ Hecke operator $T_n$ is defined by
\begin{equation}
  \label{Hecke-Eqn}
  T_n(x) := \frac1n\sum_{\alpha\in A}\psi_\alpha(x).
\end{equation}
Here $A$ is a system of representatives of those conjugacy classes of
pairs of commuting elements\footnote{In this context, ``conjugacy
  class'' refers to simultaneous conjugation.}
$(\sigma,\tau)$ of $\Sn$, with the property that the subgroup of
$\Sn$ generated by $\sigma$ and $\tau$ acts transitively on
$\{1,\dots, n\}$. For such a representative $\alpha\in A$, the
operation $\psi_\alpha$ is defined by
$$
  \psi_\alpha(x) := \on{eval}_\alpha\(\res{}{\delta}P_n(x)\),
$$
where $\delta$ denotes the inclusion of the groupoid
$X\rtimes(G\times\Sn)$
into $(X\rtimes G)\wr\Sn$, which is induced by the diagonal inclusion
of $X\rtimes G$ into $(X\rtimes G)^n$. 
(On the source of $\delta$, $\Sn$ acts
trivially on the $G$-space $X$.)   
\end{Def}

The rest of this section is dedicated to the computation of $T_n$.
We start by choosing a system of representatives $A$ as in
(\ref{Hecke-Eqn}). In order for 
$(\sigma,\tau)$ to act transitively, any two cycles of $\sigma$ must
have equal length. Let $c_1$ be a cycle of $\sigma$, w.l.o.g.
say $c_1=(1,\dots,k)$. Let $c_2$ be the $k$-cycle of $\sigma$ starting
with $\tau(1)$. Noting that $c_2$ cannot equal $c_1$, because of the
transitiveness of the action, we may assume
$c_2=(k+1,\dots,2k)$. Repeating this argument, we arrive at $N$
$k$-cycles $c_1,\dots,c_N$ of $\sigma$, where $n=Nk$, and for $1\leq
j<N$, $\tau$ 
maps the first element of $c_j$ to the first element of $c_{j+1}$. For
definiteness, we let $c_j=((j-1)k+1,\dots,jk)$. Now $\tau$ has to
map the first element of $c_N$ to an element of $c_1$, and each of the
elements of $c_1$ is possible. The conjugacy class of $(\sigma,\tau)$
is uniquely determined by $k$ and $1+m:=\tau((N-1)k+1)\in\{1,\dots,k\}$.

Fix $N$ and $k$ such that $Nk=n$, and let $\sigma$ be the element of
$\Sn$ with $N$ $k$-cycles as above. Let $g\in G$ and
$\ug:=\delta(g)=(g,\dots, g)\in G^n$. Then
$$
  P_{[\ug,\sigma]}^\str(x) = f^*_{(\ug,\sigma)}P_N^\top({}_kx\at{X^{g^k}}).
$$ 
Consider $(h,\tau)\in C_G(g)\times C_\Sn(\sigma)$,
with $\tau$ as above. In order to determine its action on
$(X^{g^k})^N$, we note that we can simplify the discusion in the proof
of Theorem \ref{Ls-Thm} by choosing all the representatives
$r_\alpha$ to equal $g^k$. Then $(h,\tau)$ acts as
$$
  (h,\dots,h, (hg^{-m},m);C)\in C_g^k\wr\Sigma_N
$$
on $(X^{g^k})^N$, where $C=(1,\dots,N)$.

Assume first that $m=0$. Then $P_N^\top$ composed with
$\on{Trace}(-,\tau)$ is the $N^{th}$ Adams operator $\psi_N$.
The set of $\tau$-fixed points of $(X^{g^k})^N$ is the image of the
diagonal map
$$
  \map{\delta'}{X^{g^k}}{(X^{g^k})^N}.
$$
Hence the $[g]$-component of $\psi_{\sigma,\tau}(x)$ is 
\begin{equation}\label{pstring-adams-Eqn}
  \(\psi_{\sigma,\tau}(x)\){[g]} = i^*(\delta')^*(\psi_N({}_kx\at{X^{g^k}})),
\end{equation}
where $$\map i{X^g\acted C_g}{X^{g^k}\acted C_{g^k}^k}$$ is the inclusion.
On coefficients, we get the following formula for characters:
$$
  \(P_n^\str(x)\)(\ug,\sigma;\uh,\tau) = \(\sum_{j\geq 0} V_j
  q^{\frac{Nj}{kl}}(g^k, h^N)\),   
$$
where $l$ is the order of $g$, and 
$x=\sum_{j\geq 0} V_j q^{\frac{j}{l}}$.
Let now $m$ be arbitrary. Then $(h,m)\in C_{g^k}^k$ acts as
$h\cdot\zeta_{lk}^{mj}$ on $V_j$. Non-equivariantly,
(\ref{pstring-adams-Eqn}) is still true, however, now the action of $h$ on
$\psi_N({}_kx\at{X^{g^k}})$ is twisted with $g^{-m}\zeta_{lk}^{mj}$.
On coefficients, we get
$$
  \(P_n^\str(x)\)(\ug,\sigma;\uh,\tau) = \(\sum_{j\geq 0} V_j
  q^{\frac{Nj}{kl}}\zeta_{kl}^{jm}(g^k,g^{-m} h^N)\),   
$$
where $x$ and $l$ are as above.
\begin{Rem}
  Viewing $x(g,h)$ as the $q$ expansion of a function $x(g,h;z)$,
  where $z$ is in the upper half-plane, and $q=e^{2\pi i z}$, we get
  (on coefficients) 
  $$
    \(T_n(x)\)(g,h;z) = \frac1n\sum_{\stackrel{Nk=n}{0\leq
        m<k}}x\(g^k,g^{-m}h^N,\frac{Nz+m}{k}\). 
  $$ 
  This is exactly the formula for the twisted Hecke operators of
  generalized Moonshine in \cite{Ganter:twistedhecke}. It follows from
  the way the twisted Hecke operators are defined in [ibid.] (namely
  using isogenies of elliptic curves)
  that our
  $P_n^\str$ are elliptic. 
  Let $G$ be the trivial group. 
  Using the theory of isogenies on the Tate-curve, 
  Matthew Ando has defined power
  operations on $K_\Tate$ in \cite{Ando:poweroperations}.
  In this case, our definition specializes to his. 
  For $G=1$ and $X$ a point, the Hecke operators
  are the usual ones acting on the $q$-expansions of modular forms.
\end{Rem}
As we noted above, our geometric picture only sees one of the two
circles of the elliptic curve. The case where $k=n$ and $N=1$
corresponds to replacing this circle by one of length $n$ and pulling
back along an $n$-fold covering map from this long circle to the short
one.
The analogous power operation for the circle we cannot see ($k=1$ adn
$N=n$), is the $n^{th}$ Adams operation.
%
%
%
%
\section{The DMVV-formula,
  Borcherds products and replicability}
We recall from \cite{Ganter:thesis} that in a cohomology theory with
power operations and a level $2$ Hopkins-Kuhn-Ravenel theory, 
the $n^{th}$ (stringy) symmetric power is defined by
$$
  \symsn(x) = \frac1{n!}\sum_{\alpha}\psi_\alpha(x),
$$
where $\alpha$ runs over all pairs of commuting elements of $\Sn$.
A similar argument to that of Section \ref{Integrality-Sec} shows that
the symmteric powers take values in $K_{\Dev,G}(X)$.
The total (stringy) symmetric power is defined as
$$\Sym_t^{\str}(x):=\sum\symsn(x)t^n.$$ 

With these definitions, the generating function argument in
\cite[9.2]{Ganter:thesis} goes through and yields:

\begin{Prop}\label{generating-function-Prop}
On the level of cohomology operations, we have
$$
  \Sym^{\str}_t(x) = \exp\(\sum_{m\geq 1}T_m(x)t^m\).
$$
\end{Prop}
In terms of elliptic curves,
this generating function can be interpreted as follows: the Hecke
operators average over the pullback along a system of representatives
for isogenies onto the elliptic curve. The $n^{th}$ symmetric power
averages over all possible $n$-fold coverings of the elliptic curve.
Each covering is made up from its
connected components, and those are isogenies. For a
detailed discussion of this picture, we refer the reader to
\cite{Ganter:twistedhecke}. 

Together with Theorem \ref{H_oo-Thm}, and applied to the case where $X$
is a point, Proposition \ref{generating-function-Prop} 
becomes
the Dijkgraaf-Moore-Verlinde-Verlinde formula:
\begin{Cor}
  Let $M$ be a compact closed complex manifold. Then
  $$
    \sum_{n\geq 0}\td^\str_\orb(M^n)t^n = \exp\(\sum_{m\geq 1}T_m(\td^\str(M))t^m\).
  $$
\end{Cor}
The right-hand side of the DMVV-formula can be rewritten as a
Borcherds product:
$$
  \prod_{i,j}\(\frac1{1-q^it^j}\)^{c(ij)},
$$
where $\phi(M)=\sum c(j)q^j$.
The best known example of a product formula is the case where the
$c(i)$ are the coefficients of the $q$-expansion of the modular function
$j-744$. 
In this case, the product is equal to the inverse of 
$$
  t(j(t)-j(q)).
$$
In other words, all the mixed terms in $t$ and $q$ in the product
$t\inv\prod(1-q^it^j)^{c(ij)}$ 
are zero.
This property encodes the replicability of the function $j-744$. Similar
identities hold for the other Moonshine functions.

Note that $j(q)-744$ equals the Witten genus of the Witten manifold
$M$ constructed in \cite{Hopkins:Mahowald}.

One rather peculiar property of these kind of product formulas is the
change of the role of $t$ from a dummy variable to a variable which
plays the same role as $q$. In our context, both the variables $t$ and
$q$ are counting winding numbers, which makes it less surprising that
they should play a similar role. Indeed, the symmetry of their roles
is clarified by
the following proposition:
\begin{Prop}
  Let $G$ be the trivial group, and view $x\in K(X)$ as an element of
  $K_\Tate(X)$. Then the total symmetric power $\Sym^\str_q(x)$ is
  Witten's exponential characteristic class
  $$
    \Sym^\str_q(x) = \bigotimes_{k\geq1}\Sym_{q^k}(x),
  $$ 
  where $\Sym_q$ stands for the total symmetric power in $K$-theory.
\end{Prop}
\begin{Pf}{}
  Note that ${}_kx=\res1{\ZZ/k}x$. Hence, for every natural number $N$,
  we have
  $$
    \langle P_N^\top({}_kx),1\rangle_{\ZZ/k\wr\Sigma_{N}} =  
    \langle P_N^\top(x),1\rangle_{\Sigma_N} = \Sym^N(x).
  $$
  The coefficient of $q^n$ in $\Sym_q^\str(x)$ is 
  \begin{eqnarray*}
    \symsn(x) &=&
    \GG{\Sn}\sum_{\sigma\tau=\tau\sigma}P_n^\str(x)\at{\delta_n}(\sigma,\tau)\\ 
    &=& \sum_{[\sigma]}\GG{C_\sigma}\sum_{\tau\in C_\sigma} \(\bigotimes
    P_{N_k}^\top({}_kx)\at{\delta_{N_k}} \)(\tau),\\
  \end{eqnarray*}
  where $[\sigma]$ corresponds to the partition $n=\sum kN_k$, so that 
  $C_\sigma\cong\prod_{k\geq 1}(\ZZ/k)\wr\Sigma_{N_k}.$
  (Since we are pulling back along the diagonal, the tensor product is
  an internal tensor product.) The expression becomes
  \begin{eqnarray*}
    \sum_{n=\sum N_k}\bigotimes_{k\geq 1}\(\GG{k^{N_k}N_k!}
    \sum_{\tau_k}\(P_{N_k}({}_k x)\at{\delta_{N_k}}\)(\tau_k)\) &=&
    \sum_{n=\sum kN_k}\bigotimes_{k\geq 1}\Sym^{N_k}(x).
  \end{eqnarray*}
  (Here $\tau_k\in\ZZ/k\wr\Sigma_{N_k}$.) This is the coefficient of
  $q^n$ in $\bigotimes_{k\geq 1}\Sym_{q^k}(x)$.
\end{Pf}
\subsection{Replicability}
Let $F(q)$ be a Laurent series with coefficients in $R(G)$ which is of
the form $$F(q)=q\inv+a_1 q+a_2 q^2+\dots.$$ From the Moonshine literature,
such $F$ are known as McKay-Thompson series. We recall\footnote{Compare
  e.g.\ \cite[(2.1)]{Teo} with $b=1$, $t=1/z$ and $F(q) = g(z)$.}
that the Faber polynomials $\Phi_{n,F}$ of $F$ are defined by
$$
  -\sum_{n=1}^\infty \Phi_{n,F} (w) t^n  = \log\(t(F(t)-w)\).
$$
Hence $\Phi_{n,F}$ is a polynomial in $w$ of degree $n$, which depends
on the first $n$ coefficients of $F$ and is uniquely characterized by
the fact that 
it is of the form 
$$\Phi_{n,F}(F(q))=q^{-n}+b_1q+b_2q^2+\dots.$$
Viewing $F(q)$ as an element in (the $[1]$-component of)
$q\inv K_{\Dev,G}(\pt)$, we arrive at the following definition:
\begin{Def}
  Let $F$ be a McKay-Thompson series. We write $F^{(a)}$ for the $a^{th}$
  Adams operator applied to $F$. We call $F$ {\em replicable}, if for
  every natural number $n$, we have
  $$
    \Phi_{n,F}(F(q)) = \sum_{\stackrel{ad = n}{0\leq b<d}}
    F^{(a)}\(\frac{a\tau +b}{d}\).
  $$  
  Here $q=e^{2\pi i\tau}$.
\end{Def}
This appears to be the right notion of replicability of McKay-Thompson
series, it is the one that turns up in \cite{Borcherds}.
Note that the right-hand side of the equation in the definition equals
$n\cdot T_n(F(q))$, where $T_n$ is the Hecke operator computed in
Section \ref{Hecke-Sec}. 
It follows immediately from the definitions and from
Corollary \ref{generating-function-Prop} that a McKay-Thompson series $F$ is
replicable if and only if it satisfies the following identity:
\begin{eqnarray*}
  F(t)-F(q) &=& t\inv\cdot\Sym^\str_t(-F(q))\\
            &=& t\inv\cdot\Lambda^\str_{-t}(F(q)),
\end{eqnarray*}
where $\Lambda^\str_{-t}(x)$ is defined as the multiplicative inverse of
$\Sym_t^\str(x)$. 
For $F(q)$ as above, $\Lambda^\str_{-t}(F(q))$ can be written as
$$
  \Lambda^\str_{-t}(F(q)) = \Lambda_{-1}\(\sum a_{m\cdot n}q^mt^n\).
$$
This is the form in which it appears in \cite[p.410]{Borcherds}.

There is a lot more to be said about the connection between power
operations in elliptic cohomology and the notion of replicability in
(generalized) 
Moonshine. We will come back 
to these topics at a 
different occasion.
\appendix
\section{Thom classes and the fixed point formula} 
\subsection{Thom classes in equivariant cohomology theories}
\begin{Def}\label{Thom-Def}
  We say that a compatible family of equivariant cohomology theories
  $\{E_G\}$ has natural Thom classes for complex vector bundles, if
  for every complex $G$-vector bundle $V$ over a pointed $G$-space $X$
  there exists a class $u_G^E(V)\in \tilde E^0(X^v)$ with the
  following properties:
  \begin{enumerate}
    \item Naturality: If $\map fXY$ is a pointed $G$-map, then
      $$
        u_G^E(f^*V) = f^*u^E_G(V).
      $$ 
    \item Multiplicativity: The family $\{E_G\}$ has external
      products, and
      if $V$ is a complex $G$-vector bundle over
      $X$ and $W$ is a complex $H$-vector bundle over $Y$, then
      $$
        u^E_{G\times H}(V\oplus W) = u^E_G(V)\boxtimes u^E_H(W).
      $$
    \item Periodicity: If $V$ is a complex $G$-representation, then
      $$
        u^E_G(V)\in \tilde E^0(\mathbb S^V)
      $$
      is a unit in $E^*$.
    \item Change of groups: If $\map \alpha HG$ is a map of groups, then
      $$
        \res{}\alpha u^E_G(V) = u_H^E(V).
      $$
  \end{enumerate}
  If $G$ is the trivial group we omit it from the notation. We might
  also omit $E$ from the notation if it is clear from the context
  which cohomology theory is meant.
\end{Def}
Fix a group $G$.
\begin{Def}\label{G-Thom-Def}
  We say that $E_G$ has natural Thom classes for complex vector
  bundles if it satisfies axioms (1) and (3) of the previous
  definition and
  $$
    (2') \quad\quad\quad u_G(V\oplus W) = u_G(V)\tensor u_G(W). 
  $$
\end{Def}
Note that (2) and (4) together imply (2').
In the situation of Definition \ref{G-Thom-Def}, let
\begin{eqnarray*}\label{alpha}
  \alpha := u(\CC)\in \tilde E^{-2}(\pt)  
\end{eqnarray*}
be the periodicity element, and 
$$
  \tau_G(V) := u_G(V)\cdot\alpha^{-d}\in\tilde E^{2d}(X^V),
$$
where $d=\dim_\CC V$.
Then the $\tau_G$ satisfy tom Dieck's axioms for equivariant Thom classes
(c.f.\ \cite[p.335]{May:Alaska}).
On the other
hand, if $E_G$ satisfies the axioms of \cite[p.335]{May:Alaska} and $E^2$
contains a unit, then the axioms of Definition \ref{G-Thom-Def} follow.
\begin{Exa}\label{MUP-def-Exa}
  Complex equivariant $K$-theory and Borel equivariant $E$ theory for even
  periodic $E$ satisfy the axioms of Definition \ref{Thom-Def}. For
  any $E$ satisfying the axioms of \cite[p.335]{May:Alaska}, 
  $$EP_G:=\bigvee_{n\in\ZZ}\Sigma^{2n} E_G$$
  satisfies the axioms of Definition \ref{G-Thom-Def}.
\end{Exa}
Recall (c.f.\ \cite{Okonek}, \cite{May:Alaska}) that $\MU_G$ is universal
among the $G$-equivariant cohomology theories with Thom classes in the
sense of tom Dieck.
\begin{Prop}\label{Okonek-Prop}
  The conditions of Definition \ref{G-Thom-Def} are equivalent to the
  existence of a map of $G$-ring spectra
  $$
    \map{\phi_G}{\MUP_G}{E_G}.
  $$
  It is the unique map of ring spectra taking Thom classes to Thom classes.
\end{Prop}
\begin{Pf}{}
  If $\phi_G$ exists, we use it to push-forward the Thom classes of
  $\MUP_G$ to $E_G$. In the situation of Definition \ref{G-Thom-Def},
  the $\tau_G$'s give rise to a map $\map{\psi_G}{\MU_G}{E_G}$, and we
  set
  $$
    \phi_G\at{\Sigma^{2n}\MU_G}:=\(\Sigma^{2n}\psi_G\)\cdot\alpha^{n}.
  $$
  Then $\phi_G$ is a map of $G$-ring spectra, and the two
  constructions are inverse to each other.
\end{Pf}
Traditionally, topologists like to work with graded rings and define
the genus corresponding to $\psi_G$ as the composite
$$
  \mathcal N_{U,G}^* \to \MU_G^*(\pt)\to E^*_G(\pt)
$$
$$
  [M]\longmapsto\psi_G(M),
$$
where the first map is the Pontrjagin-Thom map from the complex
equivariant cobordism ring to the graded coefficient ring of $\MU_G$.
Alternatively, we can consider the composite of ungraded maps
$$
  \mathcal N_{U,G}^* \to \MUP_G^0(\pt)\to E^0_G(\pt)
$$
to define $\phi_G(M)$. Hence if $[M]\in \mathcal N_{U,G}^{2d}$, then
$$
  \psi_G(M) = \phi_G(M)\alpha^{-d}.
$$
\begin{Exa}
  For the Atiyah-Bott-Shapiro $K$-theory Thom classes for complex vector
  bundles, $u_G^\td$, the periodicity element $\alpha$ is the Bott
  element $\beta$, and we have $\phi_G(M) = \td(M)$ and
  $\psi_G(M)=\td_G(M)\cdot\beta^{-d}$.   
\end{Exa}
We write $\theta_E$ for
the Thom isomorphism
$$\theta_E(x)= x\cdot u_E^G.$$
Further, $\map z{X_+}{X^V}$ denotes the zero section, and
$$
  e_G^E(V)=z^*\(u_G^E(V)\)\in E_G^0(X)
$$ 
denotes the Euler class
of $V$. Using the fact that Thom classes of trivial bundles are given
by units, one extends the definitions of the $u^E_G$ and $\theta_E$ to
virtual bundles in the usual way. Recall further that the push-forward
along the map $\map\pi M\pt$, where $M$ is a stably almost complex
oriented $G$-manifold, is defined as the composite
\begin{equation}
  \label{eq:pushforward-def}
  \pi_!\negmedspace : E_G^0(M)\to\tilde E^0_G(M^{-TM})\to\tilde
  E^0_G(\mathbb S^0),  
\end{equation}
where $TM$ is the tangent bundle of $M$, the first map is the Thom
isomorphism, and the second map is the Pontrjagin-Thom collapse. 

The genus $\phi_G$ is then computed as the push-forward of $1$:
\begin{equation}
  \label{push-forward-Eqn}
  \phi_G(M) = \pi_{G!}(1) \in\tilde E^0_G(\mathbb S^0).
\end{equation}
Similarly, the renormalized genus $\psi_G(M)$ is the graded push-forward of
$1$ which is obtained by using the graded Thom isomorphism associated
to $\tau$:
$$
  \psi_G(M) = \pi^\tau_{G!}(1)\in\tilde E^{-2d}(\mathbb S^0)
$$
(c.f.\ \cite{Okonek}, \cite{Conner:Floyd}).
\subsubsection{The $\hat A$-genus and the Todd genus Thom classes}
For the remainder of this appendix, we will closely follow the
exposition in \cite{Haynes:fpt}.
Consider the inclusion
$$
  \map{i}{\Spin(n)}{\Spinc(n)\cong\(\Spin(n)\times U(1)\)/(\ZZ/2)},
$$
and 
let $V$ be a real, $8k$-dimensional vector
bundle over $X$ which has a (chosen) $\Spin(8k)$-structure. Let 
$$
  u_{ABS}(V)\in KO(X^V)
$$
be the Atiyah-Bott-Shapiro Thom class of $V$. Forgetting the
$\Spin$-structure and viewing $V$ as a $\Spinc(8k)$-bundle via the
forgetful map $i$, we get a $K$-theory Thom class
$$u_{ABS}^\CC(V)\in K(X^V),$$
and a close look at the definition of these Thom classes on p.31 of
\cite{Atiyah:Bott:Shapiro} shows:
\begin{eqnarray*}
  u_{ABS}^\CC(V) 
                & = &  u_{ABS}(V)\tensor_\RR\CC.                
\end{eqnarray*}
We further recall the map
$$
  \map{\tilde l}{U(k)}{\Spinc(2k)},
$$
which allows us to view any complex vector bundle $V$ as a $\Spinc$
bundle. Atiyah, Bott and Shapiro prove \cite[Thm.11.6]{Atiyah:Bott:Shapiro}
$$
  u_{\td}(V) = u_{ABS}^\CC(V).
$$
These statements are summarized in the following commuting diagram of
maps of ring spectra:
\begin{equation*}
  \xymatrix{
  \on{MSU}\ar[r]\ar[d]&  \on{MSpin}\ar[r]^{\hat A}\ar[d] &
  KO\ar[d]^{\tensor\CC}\\ 
  \on{MU}\ar[r]\ar@/_4ex/[0,2]_{\td}&  \on{MSpin^\CC}\ar[r] & K \\ 
  }
\end{equation*}
Let $U^2(n)$ be the pull-back in the diagram
\begin{equation}\label{U2-Eqn}
  \xymatrix{
    U^2(n)\ar[r]^{j}\ar[d]_{\pi} & \Spin(2n) \ar[d]\\
    U(n)\ar[r] & \SO(2n). \\
  }
\end{equation}
In particular, $U^2(1)$ is the non-trivial double cover of $U(1)$ and hence
isomorphic to it.
Then $U^2(n)$ is also the pull-back in the square
\begin{equation*}
  \xymatrix{
    U^2(n)\ar[r]\ar[d]_\pi & U^2(1) \ar[d]\\
    U(n)\ar[r]^{\det} & U(1) \\
  }
\end{equation*}
(c.f.\ \cite{Haynes:fpt}).
Vector bundles with a $U^2$-structure are bundles which
have at the same time a complex structure and a Spin structure. The
second pull-back square implies that these are exactly those complex
bundles whose determinant bundle possess a square root, together with a
choice of this square root.

Note that (\ref{U2-Eqn}) gives rise to two different factorizations of
the map
$$
  U^2(n)\longrightarrow \SO(2n)
$$
through $\Spinc(2n)$:
The composition of $j$ with the inclusion $i$ becomes
$$
  \map{(j,1)}{U^2(n)}{(\Spin(2n)\times U(1))/(\ZZ/2)}.
$$
while composing $\pi$ with $\tilde l$ yields the map
$$
  \map{(j,\widetilde\det)}{U^2(n)}{(\Spin(2n)\times U(1))/(\ZZ/2)}
$$
Let $M$ be a complex Clifford module. 
Then the two above maps make $M$ into $U^2$ representations, and if
we denote the representation obtained from $(j,1)$ by
$\rho_M$ then the other one becomes $\rho_M\tensor_\CC\widetilde\det$.

Let $P$ be a principal $U^2(n)$-bundle, and let 
$$V=P\times_{U^2(n)}\RR^{2n}$$ be its associated vector bundle.
Viewing $V$ as $\Spinc(2n)$ bundle via the map $(j,1)$ (i.e.\ by first
viewing it as a $\Spin(2n)$ bundle and then forgetting some of the
structure), gives the complex Atiyah-Bott-Shapiro Thom class
$$
  u_{\hat A}^\CC(V) := u_{ABS}^\CC(V) = \chi^\CC_V\(P\times_{U^2(n)}\rho_{\mu^\CC}\).
$$
Here $\chi_V^\CC$ and $\mu^\CC$ are as defined in \cite{Atiyah:Bott:Shapiro} 
Note that this is really an abuse of notation, since $u^\CC_{ABS}(V)$
depends on the $\Spinc$ structure of $V$, not just on $V$. Viewing $V$ as
$\Spinc$-bundle via the other map $(j,\widetilde\det)$, we obtain a
different class for $u_{ABS}^\CC(V)$, namely
\begin{equation}
  \label{A-hat-Eqn}
  u_{\td}(V) = u_{ABS}^\CC(V) =\chi_V^\CC\(
  P\times_{U^2(n)}\(\rho_{\mu^\CC}\tensor_\CC\widetilde\det\)\) =
  u_{\hat A}(V)\tensor_\CC\sqrt{\det(V)}.   
\end{equation}
(The associated fibers of these two bundles are the same as complex
Clifford modules, but are viewed as $U^2(n)$-representations in
different ways. Hence the associated bundles are not necessarily
isomorphic.)  
The class $u_\td(V)$ is the Todd genus Thom class of $V$ viewed as a
complex vector bundle. If $n=4k$, we have 
$$u^\CC_{\hat A}(V)=u_{ABS}(V)\tensor_\RR\CC,$$
where $u_{ABS}(V)$ is the $KO$-Thom class of $V$ viewed as
$\Spin(8k)$-bundle. 
\subsection{Atiyah and Segal's Lefschetz fixed point formula}
Let $E$ be an even periodic cohomology theory. Fix a group $G$ and an
element $g\in G$. Then
$$
  X\mapsto E^*(X^g) =: F_G^*(X)
$$
is a $G$-equivariant cohomology theory.
Let $X$ be a $G$-pace, and 
let $V$ be a $G$-equivariant complex vector bundle over $X$. Then
the $g$-fixed points of $V$ form a vector bundle over
$X^g$, namely
$V_{g=1}$,
the eigenbundle of the eigenvalue $1$ of the action of $g$ on
$V\at{X^g}$.

Assume that we have chosen a complex orientation of $E$. Then $F_G$
inherits natural Thom classes
$$
  u^F_G(V) := u^E(V_{g=1})\in\tilde E^0((X^V)^g).
$$
Let now $E_G$ be an equivariant version of $E$ with natural Thom
classes $u^E_G$ continuing those of the complex orientation of $E$,
and assume that $E^0_\gg$ is flat over $E^0$.
Consider the natural transformation
$$
  \map r{E_G(X)}{E_\gg(X^g)}
$$
defined by
$$
  r:=i^*\circ\res G\gg,
$$
where $\map i{X^g}X$ is the inclusion of the fixed points.
\begin{Def}
  For a complex $G$-vector bundle $V$ over $X$, we write $V_{g=\zeta}$ for the
  $\zeta$-eigenbundle of the action of $g$ on $V\at{X^g}$. Write
  $V_{g\neq 1}$ for the orthogonal complement of $V_1$,
  $$
    V_{g\neq 1} = \bigoplus_{\zeta\neq 1}V_\zeta.
  $$
\end{Def}
\begin{Exa}
  Let $M$ be a smooth complex $G$-manifold, and let $T$ its
  tangent bundle. Then 
  $T_{g=1}$ is the tangent bundle $TM^g$ of $M^g$, and $T_{g\neq
  1}$ is the normal bundle $N^g$ of $i$.    
\end{Exa}
\begin{Prop}\label{Lefschetz-Prop}
  For a complex $G$-vector bundle $V$ over $X$, we have
  $$
    r\(\theta_E(x)\) = \theta_F(r(x)\cdot e_G^E(\Vrest)). 
  $$
\end{Prop}
Its proof relies on the following lemma about relative zero
sections: 
\begin{Lem}[{\cite{Rudyak}}]\label{eu-Lem}
  Let $V$ and $W$ be (eqivariant) vector bundles over $X$, let $$\map
  sWV\oplus W$$ be the inclusion, and write $X^s$ for the induced map
  of Thom spaces $$\map{X^s}{X^W}{X^{V\oplus W}}.$$
  Then $X^s$ pulls back $u_G(V\oplus W)$ to $$e_G(V)\cdot u_G(W)\in\tilde
  E^0_G(X^W).$$ 
\end{Lem}

\begin{Pf}{of Proposition \ref{Lefschetz-Prop}}
  The fixed points
  inclusion $$\map i{(X^V)^g}{X^V}$$ factors into the composition of
  the two maps
  $$\map{i_1:=(X^g)^s}{(X^g)^{V_{g=1}}}{(X^g)^{\Vrest\oplus V_{g=1}}}$$
  and
  $$\map{i_2}{(X^g)^{\Vrest\oplus V_{g=1}}}{X^V},$$
  where $i_2$ is the map of Thom spaces obtained from the fixed points
  inclusion $\map{i_g}{X^g}{X}$ and the corresponding bundle map 
  $$
    \map{\tilde i}{i^*(V)\cong \Vrest\oplus V_{g=1}}{V}.
  $$
  We have
  \begin{eqnarray*}
    r(\theta_E(x)) & = & r\(x\cdot u_G^E(V)\)\\
    & = & r(x)\cdot r\(u_G^E(V)\)\\
    & = & r(x)\cdot i_1^*\circ i_2^*\res{G}{\gg}u_G^E(V)\\
    & = & r(x)\cdot i_1^*\circ i_2^* u_\gg^E(V)\\
    & = & r(x)\cdot i_1^* \(u_\gg^E(\Vrest\oplus V_{g=1})\)\\
    & = & r(x)\cdot e_\gg^E(\Vrest)\cdot u_\gg^E(V_{g=1})\\
    & = & \theta_F(r(x)\cdot e_\gg^E(\Vrest)),
  \end{eqnarray*}
  where the second to last equality is Lemma \ref{eu-Lem}.
  Let $\alpha$ be the unique map from $\gg$ to the trivial group.
  Since the action of $\gg$ on $V_{g=1}$ is trivial, we have
  $$u_\gg^E(V_{g=1}) =\res{}\alpha u^E(V_{g=1})=u_G^F(V).$$
\end{Pf}
The correction factor $e_\gg(\Vrest)$ is an exponential characteristic
class. If we assume the Euler classes of $\gg$-representations to be
invertible and set 
$$
  e_\gg(-V):= e_\gg(V)\inv\in\tilde E^0_\gg(\mathbb S^V),
$$ 
it follows that Proposition \ref{Lefschetz-Prop} holds for virtual
bundles, too. 
In the case $E_G = K_G$, and $u^E_G = u^\td_G$,
the Riemann-Roch formula in \cite{Dyer} 
$$
  \pi_!^E(x) = \pi_!^F\(r(x)\cdot e^\td_\gg(-N^g)\)
$$

becomes exactly the
Atiyah-Segal result \cite[2.10]{Atiyah:Segal}
\begin{equation}
  \label{Lefschetz-Eqn}
  \(\on{ind}_G^X(x)\)_g = \(\on{ind}^{X^g}_1\tensor\id_{R(\gg)}\)_g
  \(\frac{r(x)}{\lambda_{-1}(N^g)}\).  
\end{equation}
Here $\lambda_{-1}(V)$ denotes the alternating sum of the exterior
powers of the $\gg$-representation $V$.
Let now $u^\AA$ be the complex Atiyah-Bott-Shapiro Thom class for
$\Spin(2d)$-vector bundles. For simplicity, we assume that $M$ is a
$U^2$-manifold. 
We combine the classical
Riemann-Roch formula for the Chern character
\begin{eqnarray*}
  \pi_!^\AA(a) & = &
  \pi^H_!\left(\frac{\on{ch}(e^\AA(-TM))} {e^H(-TM)}\cdot\on{ch}(a)\right) \\
  &=& \int_M \widehat{\mathcal A}(TM)\on{ch}(a)
\end{eqnarray*}
with the Lefschetz formula (\ref{Lefschetz-Eqn}) to obtain
\begin{eqnarray}
  \label{A-hat-Lefschetz-Eqn}
  \notag 
  \({}^M\pi_{G!}^\AA(a) \)(g) &=& 
  \int_{M^g}\frac{e^H(TM^g)}{\on{ch}(e^\AA(TM\at{M^g})(g))}
  \cdot\on{ch}((a\at{M^g})(g))\\
  &=&
  \int_{M^g}\frac{e^H(TM^g) \on{ch}\(\sqrt{\det TM\at{M^g}}(g)\)}
  {\on{ch}(e^\td(TM\at{M^g})(g))}
  \cdot\on{ch}((a\at{M^g})(g))\\
  \notag &=&
  \varepsilon\(\(\prod_{x_i}x_i\)\(\prod_{r=0}^{k-1}\prod_{y_i}
  \frac{e^{\frac{y_i}2+\pi i\frac rk}}{1-e^{y_i+2\pi i\frac
      rk}}\)\cdot\on{ch}((a\at{M^g})(g))\)\left[ M^g\right]. 
\end{eqnarray}
Here the $x_i$ are the Chern roots of $TM^g$, for fixed $r$, the $y_i$
run over the 
Chern roots of the $e^{2\pi i\frac rk}$ eigenbundle $\(TM\at{M^g}\)$,
and $\varepsilon =\pm 1$ depends on the way that the action of $g$ on
$\det TM$ lifts to an action on $\sqrt{\det TM}$ and on our choice of
identification of $U^2(1)$ with $U(1)$ above. 
\begin{Rem}
  The above discussion remains valid if we replace $\gg$ by a
  topologically cyclic group e.g.\ $\SS$ or $\Rl$.
\end{Rem}
Let now $P$ be a principal $\Spin(2n)$ bundle over $X$, and
$V=P\times_{\Spin(2n)}\RR^{2n}$ its associated vector bundle.
We assume that $V$ has an even $\SS$-action, i.e., one that is
induced by an $\SS$-action on $P$.
We write 
$$
  u^{\hat A}_\SS(V) := u_\SS^{\hat A,\CC}(V) 
$$
for its equivariant \ABS Thom class in coplex $K$-theory.
Note that $V_{g=1}$ and its orthogonal complement are still even
dimensional $\Spin$ bundles. Hence the above discussion goes through
with $U(n)$ replaced by $\Spin(2n)$, and if we let $\pi_!^{\hat A}$
denote the (equivariant) push-forward in $K$-theory defined by 
$u^{\hat A}$, we get
$$
  {}^{\hat A}\pi_{\SS!}^{X}(x) = {}^{\hat A}\pi_!^{X^\SS}\(r(x)\cdot
  e_\SS^{\hat A}(-N)\),
$$ 
where $N$ is the normal bundle of the fixed point inclusion.
In order to compute the equivariant Euler class of $N$ (or more
generally of $V_{g\neq 1}$), recall that, since
$\SS$ acts fibre preserving and fixed-point free on $N$, we have by
(\ref{decomp-Eqn})
\begin{equation}
  \label{normal-Eqn}
  N\cong_\RR\bigoplus_{j\geq 1}V_j q^j,  
\end{equation}
and $N$ can be equipped (uniquely) with a complex structure making
this a complex 
isomorphism. 
Then (\ref{A-hat-Eqn}) gives
\begin{eqnarray}
  \notag
   e^{\AA}_\SS(N) & = &e_\SS^\td(N)\cdot\(\sqrt{\det(N)}\acted\SS\)\inv\\
  \label{Todd-A-hat-Eqn}
   &=& q^{-{d}} \(\sqrt{\det N}\)\inv \prod_{j\geq 1}\Lambda_{-q^j}V_j,
\end{eqnarray}
where $d=\sum_j j\dim_\RR V_j $.

\bibliographystyle{alpha}
\bibliography{stringy.bib}
\end{document}

%% file: moonshine.macros.tex
\theoremstyle{plain}
\newtheorem {Thm} {Theorem}[section]
\newtheorem* {Thm*} {Theorem}
\newtheorem* {Prop*} {Proposition}
\newtheorem {Lem}[Thm] {Lemma}
\newtheorem {Prop}[Thm] {Proposition}
\newtheorem {Cor}[Thm] {Corollary}

\theoremstyle {definition}
\newtheorem {Def}[Thm] {Definition}

\newtheorem {Rem}[Thm] {Remark}
\newtheorem {Exa}[Thm] {Example}
\newenvironment{Pf}[1]{{\noindent\sc Proof #1:}}{\qed\\}
\newcommand {\specialmap} [4] {\text {$ #1\negmedspace : #2 #3 #4 $}}
\newcommand {\map} [3] {\specialmap {#1} {#2}{\to} {#3}}

\newcommand {\id} {\operatorname{id}}

\newcommand {\at}[1] {\arrowvert_{#1}}

\newcommand {\ind} {\operatorname{ind}}
\renewcommand {\(} {\left(}
\renewcommand {\)} {\right)} 

\renewcommand {\geq} {\geqslant}
\newcommand {\Oplus} {{\bigoplus\limits_{n\geq 0}}}
\newcommand {\cOplus} {\hat\Oplus}

\def\acted{\hspace{.1cm}{
        \setlength{\unitlength}{.32mm}
        \linethickness{.09mm}
        \begin{picture}(8,8)(0,0)
        \qbezier(1,6)(3.5,8.3)(6,7)
        \qbezier(6,7)(9.5,4)(6,1)
        \qbezier(6,1)(3.5,-.3)(1,2)
        \qbezier(1,6)(1.9,7.5)(1.2,9)
        \qbezier(1,6)(3,6.1)(3.8,4.4)
        \end{picture}\hspace{.1cm}
        }}

\newcommand {\sub} {\subseteq}

\newcommand {\CC} {\mathbb C}
\newcommand {\on}[1] {\operatorname{#1}}

\newcommand {\EG}[1] {{E_G(#1)}}

\newcommand {\eps} {\varepsilon}

\newcommand {\forb}[1] {{\phi_{\on{orb}}(#1)}}
\renewcommand {\gg} {{(g_1,\dots,g_h)}}
\newcommand {\GG}[1] {{\frac{1}{|#1|}}}

\newcommand {\GS} {{G\wr\Sn}}

\newcommand {\hinf} {{H_\infty}}

\renewcommand {\ind}[2] {\on{ind}\arrowvert_{#1}^{#2}}

\newcommand {\Korb}[1] {K_{\on{orb}}(#1)}

\renewcommand {\leq} {\leqslant}
\newcommand {\LG} {\Lambda(G)}

\newcommand {\mmod} {/\!\!/}

\newcommand {\MU} {{\on{MU}}}

\newcommand {\orb} {{\on{orb}}}

\newcommand {\pt} {\on{pt}}
\newcommand {\Ptop} {{P^{\on{top}}}}

\newcommand {\res}[2] {\on{res}\arrowvert^{#1}_{#2}}

\newcommand {\Sj} {{\Sigma_j}}

\newcommand {\Sn} {{\Sigma_n}}
\newcommand {\Sm} {{\Sigma_m}}

\newcommand {\SO} {{\Sigma^\infty_+\Omega^\infty}}

\renewcommand {\SS} {{\mathbb{S}^0}}

\newcommand {\Stop} {{S_t^{\on{top}}}}

\newcommand {\Sym} {\on{Sym}}
\newcommand {\tensor}{\otimes}
 
\newcommand {\Tate} {{\on{Tate}}}

\newcommand {\td} {{\on{Td}}}
\renewcommand {\top} {{\on{top}}}

\newcommand {\tr}[1] {\on{Trace}(g | #1)}

\newcommand {\XX} {{X(1)}}
\newcommand {\LX} {\Lambda(\XX)}

\newcommand {\ZZ} {\mathbb Z}

\newcommand {\lps} {[\! [}
\newcommand {\rps} {]\! ]}

\newcommand {\ps}[1] {\lps #1\rps}